\documentclass[a4paper,12pt]{amsart}

\usepackage{amssymb}

\def\q{\footnote}
\def\sa{{\sf{A}}}
\def\ss{{\sf{S}}}
\def\sp{{\sf{P}}}
\def\se{{\sf{eu}}}
\def\kz{{\tt{KZ}}}
\def\reg{_{{{\operatorname{reg}\nolimits}}}}
\def\tor{_{{{\operatorname{tor}\nolimits}}}}
\def\ii{^\dag}
\def\ccirc{{{}_{\,^{^\circ}}}}

\def\Bbar{\bar{B}}

\def\Mbar{\bar{M}}

\def\Pbar{\bar{P}}
\def\Qbar{\bar{Q}}

\def\Xbar{\bar{X}}

\def\gV{{\mathbb{V}}}

\def\iso{\buildrel \sim\over\to}

\def\Gg{{\mathfrak{g}}}

\def\Gm{{\mathfrak{m}}}

\def\CA{{\mathcal{A}}}
\def\CB{{\mathcal{B}}}
\def\CC{{\mathcal{C}}}
\def\CD{{\mathcal{D}}}
\def\CE{{\mathcal{E}}}
\def\CF{{\mathcal{F}}}

\def\CH{{\mathcal{H}}}

\def\CJ{{\mathcal{J}}}

\def\CO{{\mathcal{O}}}
\def\CP{{\mathcal{P}}}

\def\CU{{\mathcal{U}}}

\def\CW{{\mathcal{W}}}

\def\BC{{\mathbf{C}}}

\def\BQ{{\mathbf{Q}}}

\def\BZ{{\mathbf{Z}}}
\def\Bk{{\mathbf{k}}}

\def\eps{\varepsilon}

\def\Ch{\operatorname{Ch}\nolimits}

\def\mcoh{\operatorname{\!-coh}\nolimits}

\def\End{\operatorname{End}\nolimits}
\def\Ext{\operatorname{Ext}\nolimits}

\def\GL{\operatorname{GL}\nolimits}

\def\gr{{\operatorname{gr}\nolimits}}

\def\Hom{\operatorname{Hom}\nolimits}

\def\Homgr{\operatorname{Homgr}\nolimits}

\def\Id{\operatorname{Id}\nolimits}

\def\Ind{\operatorname{Ind}\nolimits}

\def\Irr{\operatorname{Irr}\nolimits}

\def\mMod{\operatorname{\!-mod}\nolimits}

\def\mMOD{\operatorname{\!-Mod}\nolimits}

\def\opp{{\operatorname{opp}\nolimits}}

\def\mproj{\operatorname{\!-proj}\nolimits}

\def\Res{\operatorname{Res}\nolimits}

\def\mtilt{\operatorname{\!-tilt}\nolimits}

\def\ie{{\em i.e.}}

\def\tM{{\tilde{M}}}

\def\tDelta{{\tilde{\Delta}}}

\def\tCO{{\tilde{\CO}}}

\def\uc{\operatorname{\underline{c}}\nolimits}

\newtheorem{thm}{Theorem}[section]
\newtheorem{lemma}[thm]{Lemma}
\newtheorem{cor}[thm]{Corollary}
\newtheorem{prop}[thm]{Proposition}
\newtheorem{defi}[thm]{Definition}

\theoremstyle{definition}
\newtheorem{rem}[thm]{Remark}
\newtheorem{hyp}{Hypothesis}

\input xypic
\xyoption{all}
\def\COln{\CO^{\mathrm{ln}}}

\begin{document}
%\date{30 October 2002}
\title{On the category $\CO$ for rational Cherednik algebras}
\author{Victor Ginzburg,\, Nicolas Guay,\, Eric Opdam and\, Rapha\"el Rouquier}
\maketitle
\begin{abstract}{We study   the category $\CO$ of representations of
the rational Cherednik algebra $\sa_W$
attached to a complex reflection group $W$.
We construct an exact functor,
called  {\it Knizhnik-Zamolodchikov functor}: $\CO\to \CH_W\mMod$,
where $\CH_W$ is  the (finite) Iwahori-Hecke algebra associated to $W$.
We prove that the Knizhnik-Zamolodchikov  functor induces an
equivalence between $\CO/\CO\tor$, the quotient of $\CO$ by the subcategory
of $\sa_W$-modules supported on the discriminant,
and  the
category of finite-dimensional $\CH_W$-modules. 
The standard
$\sa_W$-modules
go, under this equivalence, to certain modules arising in
Kazhdan-Lusztig theory of ``cells'', provided
 $W$ is a Weyl group
and the Hecke
algebra 
$\CH_W$ has equal parameters.
We prove that the category $\CO$ is
equivalent to the module category over a finite dimensional algebra,
a generalized "$q$-Schur algebra" associated to $W$.}
\end{abstract}
{\footnotesize
\tableofcontents}
\section{Introduction}
Let $W$ be  a complex reflection group  acting on a vector space $V$. 
Let  $\sa_W$ denote the
rational Cherednik algebra introduced in \cite{EtGi} as
a certain deformation of $\CD(V)\rtimes W$, the cross-product of $W$
with the  algebra of polynomial differential operators on $V$.
The algebra $\sa_W$ can be also realized as an algebra of operators (Dunkl
operators)
acting on polynomial functions on $V$.
When $W$ is a Weyl group, $\sa_W$ is a rational degeneration of the double affine
Hecke algebra.

\smallskip
A nice category $\CO$ of $\sa_W$-modules has been
discovered in \cite{DuOp}, cf. also \cite{BeEtGi}. It
shares many similarities with the Bernstein-Gelfand-Gelfand
category $\CO$ for a finite-dimensional
semi-simple Lie algebra.

We develop a general approach to the category $\CO$ for a rational
Cherednik algebra, similar in spirit
to Soergel's analysis, see \cite{So1}, of the category $\CO$ in the Lie
algebra case. Specifically, in addition to the algebra $\sa_W$,
we consider an appropriate (finite) Hecke algebra $\CH_W$,
and construct an exact functor $\kz: \CO \to \CH_W\mMod$,
that may be thought of as a Cherednik algebra analogue of the
functor $\mathbb{V}$ of  \cite{So1}.
One of our main results  says that
the functor $\kz$ is fully faithful on projectives.
Thus, the (noncommutative!) Hecke algebra plays, in our case,
 the role similar to that the coinvariant algebra (= cohomology of the flag manifold)
 plays
in the Lie algebra case. It is also interesting to note that, in both cases,
the algebra in question is Frobenius.

To prove our results, in \S \ref{sectioncatO} we develop some basic
representation 
theory  over a ground ring
(which is not necessarily a field) of  a general associative algebra with a 
{\em triangular decomposition}.
 This generalizes earlier work
of the second author \cite{Gu} and
of the last two authors (unpublished).
 Such generality will be essential for us in order to use deformation
arguments in \S \ref{sectionHecke}.
The results of section \ref{sectioncatO}
are applied to Cherednik algebras in \S \ref{catOCherednik}.

\smallskip
In \S \ref{sectionDmod}, we explain how to generalize some classical
constructions for $\CD(V)$, the Weyl algebra, (such as characteristic varieties, duality)
to the rational Cherednik algebra.
We study two kinds of dualities. One of them is related to Fourier
transform while the other, much more important
 one,  generalizes the usual
(Verdier type) duality on
$\CD$-modules. This enables us to show that the
Ringel dual of category $\CO$ is a category $\CO$ for the dual
reflection group. We also give a formula for the dimension of the
characteristic variety involving only the highest weight structure of
$\CO$.

\smallskip
Our most important results are concentrated
in \S \ref{main}. We use the de Rham functor
for Knizhnik-Zamolodchikov type $\CD$-modules over the complement of
the ramification locus in $V$. This way, we
 relate the category $\CO$ with a Hecke algebra.
We prove that the category $\CO$ can be recovered from
its quotient by the subcategory of objects with non-maximal
characteristic variety (Theorem \ref{doublecentralizer}
and Corollary \ref{doublecentralizer_cor}).

Then, we obtain  a ``double centralizer" Theorem \ref{OfromH3},
asserting in particular
that the category $\CO$ is equivalent to the category of
modules over the endomorphism ring of some Hecke algebra module.
 A crucial point is the proof that the
de Rham functor sends the $\CD$-modules coming from objects of $\CO$
to representations of the braid group that factor through the Hecke algebra
(Theorem \ref{factorHecke}).\smallskip

In a different perspective, our results provide a solution to the problem of
associating a generalized ``$q$-Schur algebra'' to an arbitrary
finite complex reflection group $W$. This seems to be new even when
$W$ is a Weyl group (except for types $\bf{A}, \bf{B}$).
For instance, let $W$ be the Weyl group of an irreducible simply-laced
root system. Then, the data defining the Cherednik algebra $\sa_W$
reduces to a single complex parameter $c\in\BC$. In this case,
$\CH_W$ is the standard Iwahori-Hecke algebra of $W$, specialised at the
parameter $q=e^{2\pi i c}$. If $c$ is a rational number, then
$q$ is a root of unity, and the corresponding category
$\CH_W\mMod$ becomes quite complicated. Our results show that
the category $\CO$ for  $\sa_W$ may be viewed as a natural
``quasi-hereditary cover'' of the category $\CH_W\mMod$,
which is not itself quasi-hereditary.
As a consequence, the decomposition matrices of Hecke algebras
(in characteristic $0$) are triangular (Corollary \ref{decHecke}).
We remark that,
in view of \cite{CPS2}, one might have expected on
general grounds
that  the category $\CH_W\mMod$ only has a ``stratified cover'',
which is weaker than having a ``quasi-hereditary cover''.

The reader should be reminded that, in type  $\bf{A}$,
a well-known ``quasi-hereditary cover'' of $\CH_W\mMod$ is provided by
the $q$-Schur algebra. We expect that the latter category is
equivalent to the category $\CO$.
Furthermore, for an arbitrary finite Weyl group $W$, we prove in \S
\ref{subKL}
that the $\kz$-functor sends
 the standard modules in $\CO$
to  modules over
the Hecke algebra (with equal parameters) 
that can be described via Kazhdan-Lusztig's theory of cells.
 It follows in particular that, in type $\bf{A}$, the standard modules in category $\CO$
go to 
Specht (or `dual Specht', depending on the sign of parameter `$c$')
$\CH_W$-modules, introduced in \cite{DJ}.

\bigskip
{\bf
{Acknowledgments.}} {\footnotesize{
The second named author gratefully acknowledges the financial support of
the Fonds NATEQ. The third named author was partially supported by a
Pioner grant of 
the Netherlands Organization for Scientific Research (NWO).}}

\section{Category $\CO$}
\label{sectioncatO}
\subsection{Algebras with triangular decomposition.}\label{triag}
In this section, we assume  given  an associative algebra
$A$  with a  triangular decomposition.
We study a category $\CO(A)$ of $A$-modules,
 similar to the Bernstein-Gelfand-Gelfand category $\CO$ for
a complex semi-simple Lie algebra.
The main result of this section is Theorem
\ref{highestweightcategory} below, saying that 
the category $\CO(A)$ is a highest weight category (in the sense of
\cite{CPS1}).

Throughout this section \ref{sectioncatO}, let $k_0$ be an algebraically
closed field and  $k$ a commutative noetherian $k_0$-algebra.

Let $A$ be a graded $k$-algebra with three graded subalgebras
$B$, ${\Bbar}$ and $H$ such that
\begin{itemize}
\item $A={\Bbar}\otimes H\otimes B$ as $k$-modules
\item $B$ and ${\Bbar}$ are projective over $k$
\item $B\otimes H=H\otimes B$ and $H\otimes {\Bbar}={\Bbar}\otimes H$
\item $B=\bigoplus_{i\leq 0}\, B_i$, ${\Bbar}=\bigoplus_{i\geq 0}\,{\Bbar}_i$, and $B_0={\Bbar}_0=k$ and
$H\subset A_0$.
\item $H=k\otimes_{k_0}H(k_0)$ where $H(k_0)$ is a finite dimensional
semi-simple split $k_0$-algebra
\item
the grading on $A$ is inner, \ie, there exists $\partial\in A_0$ such
that $A_i=\{u\in A | \partial u-u\partial = iu\}$.
\end{itemize}

We denote by $BH$ and ${\Bbar}H$ the subalgebras $B\otimes H$ and ${\Bbar}\otimes H$.
We put $B^i=B_{-i}$.
We denote by $\Irr(H(k_0))$ the set of isomorphism classes of finite dimensional
simple $H(k_0)$-modules.
We put $\partial=\partial'-\partial_0$
with $\partial'\in {\Bbar}\otimes H\otimes B^{>0}$ and $\partial_0\in Z(H)$.
For $E\in\Irr(H(k_0))$, we denote by $c_E$ the scalar by which $\partial_0$
acts on $k\otimes_{k_0} E$.

\smallskip
The theory developped here is closely related to
the one developped by
Soergel \cite[\S 3-6]{So2} in the case where
$\Gg$ is a graded Lie algebra with $\Gg_0$ reductive,
$A=\CU(\Gg)$, $B=\CU(\Gg_{>0})$, ${\Bbar}=\CU(\Gg_{<0})$ and
$H=\CU(\Gg_0)$.\footnote{In the Lie algebra case, the algebra $H=\CU(\Gg_0)$
 is not finite dimensional. One then has to restrict oneself to 
the consideration 
 of $H$-semisimple $A$-modules only. 
The theory developed below easily extends to such a case.}

\pagebreak[3]\subsection{Locally nilpotent modules}

We denote by $\COln$ the full subcategory of the category of
$A$-modules consisting of those modules that are
locally nilpotent for $B$, i.e.,
an $A$-module $M$ is in $\COln$ if for every $m\in M$,
there exists  $n\gg 0$ such that $B^{> n}\cdot m=0$.
 This is a Serre subcategory of the category of $A$-modules.

\begin{rem}
The canonical functor $D^b(\COln)\to D^b(A)$ is not faithful in general.
Nevertheless,  for $i=0,1$, and any  $M, M'\in \COln$, one still has
$\Ext^i_{\COln}(M,M')\iso \Ext^i_A(M,M')$.
\end{rem} 

\pagebreak[3]\subsection{Standard modules} 
\subsubsection{}
Let $h\in H$.
We denote by $\phi_h:{\Bbar}\to {\Bbar}\otimes H\subseteq A$ the map defined by
$\phi_h(\bar{b})=h\otimes \bar{b}$. Similarly, we denote by
$\psi_h:B\to H\otimes B\subseteq A$ the map defined by
$\psi_h(b)=b\otimes h$.

Let $E$ be an $H$-module. The augmentation
$B\to B/B^{>0}=k$ induces a morphism of algebras
$BH\to H$ and we view $E$ as a $BH$-module by restriction via this morphism.
All simple $BH$-modules that are locally nilpotent over $B$
are obtained by this construction, starting with
$E$ a simple $H$-module.

We put 
$$\Delta(E)=\Ind_{BH}^A E= A\otimes_{BH} E.$$

The canonical isomorphism $\Delta(E)\iso {\Bbar}\otimes E$
is an isomorphism of graded ${\Bbar}H$-modules ($E$ is viewed in degree $0$),
where ${\Bbar}$ acts by multiplication on ${\Bbar}$ and the action of $h\in H$ is
given by $\phi_h\otimes_H 1_E:{\Bbar}\otimes E\to {\Bbar}\otimes H\otimes_H E={\Bbar}\otimes E$.

\smallskip
We now put $\nabla(E)=\Homgr_{{\Bbar}H}^\bullet(A,E)=
\bigoplus_i \Homgr_{{\Bbar}H}^i(A,E)$ (this is also the submodule of elements of
$\Hom_{{\Bbar}H}(A,E)$ that are locally finite for $B$). Here, $E$ is viewed
as a ${\Bbar}H$-module via the canonical morphism ${\Bbar}H\twoheadrightarrow
({\Bbar}/{\Bbar}_{>0})\otimes H=H$.

We have an isomorphism of graded $BH$-modules
$\nabla(E)\iso\Hom_k(B,k)\otimes E$
where $B$ acts by left multiplication on $\Hom_k(B,k)$ and the
action of $h\in H$ is given by $f\otimes e\mapsto
(b\otimes e\mapsto (1\otimes f)(\psi_h(b))e)$.

The $A$-module $\Delta(E)$ is a graded module, generated by its degree $0$
component.
The $A$-module $\nabla(E)$ is also graded. Both $\Delta(E)$
and $\nabla(E)$ are concentrated in non-negative degrees, hence
are locally nilpotent for $B$.

\subsubsection{}

We have
$$\Ext^i_A(\Delta(E),\nabla(F))\simeq
\Ext^i_{{\Bbar}H}(\Res_{{\Bbar}H}\Delta(E),F)\simeq
\Ext^i_{{\Bbar}H}(\Ind_H^{{\Bbar}H}E,F)\simeq \Ext^i_H(E,F).$$
It follows that, when $k$ is a field and $E,F$ are simple, then
\begin{equation}
\label{Extvanishing}
\Ext^i_A(\Delta(E),\nabla(F))=0\text{ if }i\not=0\text{ or }E\not\simeq F
\quad\text{ and }\quad\Hom_A(\Delta(E),\nabla(E))\simeq k.
\end{equation}

Let $N$ be any $A$-module. We have
\begin{equation}
\label{adjunctionDelta}
\Hom_A(\Delta(E),N)\iso \Hom_{BH}(E,\Res_{BH}N)
\end{equation}

\subsubsection{}
A $\Delta$-filtration for a $A$-module $M$ is a filtration
$0=M_0\subset M_1\subset\cdots\subset M_n=M$
with $M_{i+1}/M_i\simeq \Delta(k\otimes_{k_0} E_i)$ for some
$E_i\in\Irr(H(k_0))$. We denote by $\CO^\Delta$ the full subcategory
of $\COln$ of objects with a $\Delta$-filtration.

Given an $H$-module $E$ and $n\ge 0$, we also consider
 more general modules 
$$\Delta_n(E)=\Ind_{BH}^A \bigl((B/B^{> n})\otimes_k E\bigr)$$
The modules $\Delta_n(k\otimes_{k_0} F)$ have a $\Delta$-filtration,
when $F$ is a finite dimensional $H(k_0)$-module.

For $N$ a $A$-module, we have
$$\Hom_A(\Delta_n(E),N)\iso \Hom_{BH}\bigl((B/B^{> n})\otimes_k E,N\bigr).$$
As a consequence, we have a characterization of $B$-locally nilpotent
$A$-modules~:

\begin{prop}
\label{charlocnilp}
Let $N$ be a $A$-module. Then, the following are equivalent
\begin{itemize}
\item
$N$ is in $\COln$
\item
$N$ is a quotient of a (possibily infinite) sum of $\Delta_n(E)$'s
\item
$N$ has an ascending filtration whose successive quotients are quotients
of $\Delta(E)$'s.
\end{itemize}
\end{prop}

\pagebreak[3]\subsection{Graded modules}
\label{sectiongradedmodules}

\subsubsection{}
Given $\alpha\in k$ and $M$ a $A$-module, define generalized weight
spaces in $M$ by
$$\CW_\alpha(M)=\{m\in M \;|\; (\partial-\alpha)^nm=0\text{ for }n\gg 0\}.$$

Let $\CO$ be the full subcategory of $\COln$ consisting
of those modules $M$ such that
$M=\sum_{\alpha\in k} \CW_\alpha(M)$ where
$\CW_\alpha(M)$ is finitely generated over $k$, for every $\alpha\in k$.
This is a Serre subcategory of the category of $A$-modules.

\smallskip
Let $\tCO$ be the category of graded $A$-modules
that are in $\CO$.
This is a Serre subcategory of the category of graded $A$-modules.

Let $\tCO^\alpha$ be the full subcategory of $\tCO$ consisting of those
objects $M$ such that $M_i\subseteq \CW_{i-\alpha}(M)$ for all $i$.
Note that this amounts to requiring that
$\partial'-(i+c_F-\alpha)$ acts nilpotently on
$\Homgr^i_{H}(k\otimes_{k_0} F,M)$
for $F\in\Irr(H(k_0))$, since $\partial$ and $\partial_0$ commute.

More generally, if $I$ is a subset of $k$, we denote by $\tCO^I$
the full subcategory of $\tCO$ consisting of those
objects $M$ such that $M_i\subseteq \sum_{\alpha\in I}\CW_{i-\alpha}(M)$.

\smallskip
We denote by $\tDelta(E)$ the graded version of $\Delta(E)$
(it is generated in degree $0$ and has no terms in negative degrees).
Further, write $\langle r\rangle$ for `grading shift by $r$' of
a graded vector space.

\begin{lemma}
\label{actionE}
Let $E\in\Irr(H(k_0))$.
We have $\tDelta(k\otimes_{k_0} E)\langle r\rangle\in
\tCO^{c_E-r}$.
\end{lemma}

\begin{proof}
Note that $\partial'$ acts as zero on $\tDelta(k\otimes_{k_0} E)_0$,
since $B^{>0}$
acts as zero on it. So, $\partial$ acts as $-c_E$ on it.
It follows that $\partial$ acts by $i-c_E$ on
$\Bbar_i\tDelta(k\otimes_{k_0} E)_0=\tDelta(k\otimes_{k_0} E)_i$ and we are
done.
\end{proof}

\subsubsection{}
\label{sectiondecomposition}

Let $\CP$ be the quotient of $\bigcup_{E\in\Irr(H(k_0))}(c_E+\BZ)$ by the
equivalence relation given as the transitive closure of the relation~:
$\alpha\sim \beta$ if $\alpha-\beta$ is not invertible.

\smallskip
We make the following assumption until the end of \S \ref{sectiongradedmodules}.

\begin{hyp}
\label{width}
We assume that $c_E\sim c_E+n$ for some $n\in\BZ$
implies $n=0$ (this holds for example when $k$ is a local ring of
characteristic zero).
\end{hyp}

\begin{prop}
\label{decomposition}
We have $\tCO=\bigoplus_{a\in \CP}\tCO^a$.

The image by the canonical functor $\tCO\to\CO$ of
$\tCO^{a+n}$ is a full subcategory $\CO^{a+\BZ}$ independent of $n\in\BZ$.

We have $\CO=\bigoplus_{a\in \CP/\BZ}\CO^{a+\BZ}$ and
the forgetful functor $\tCO^a\iso  \CO^{a+\BZ}$ is an equivalence.
\end{prop}

\begin{proof}
Let $M$ be an object of $\CO$.
Let $a\in \CP$ and $M^a=\sum_{\alpha\in -a+\BZ}\CW_{\alpha}(M)$.
By Lemma \ref{actionE} and Proposition \ref{charlocnilp},
we have a decomposition $M=\bigoplus_{a\in \CP/\BZ} M^a$ as
$A$-modules.

Similarly, given $\tM\in\tCO$, we have
$\tM=\bigoplus_{a\in\CP} \tM^a$ where
$$\tM^a=\bigoplus_i \sum_{\alpha\in a}(\CW_{i-\alpha}(M) \cap M_i)
\in\tCO^a.$$

Given $M\in\CO^{a+\BZ}$, we put a grading on $M$ by setting
$M_i=\sum_{\alpha\in i-a}\CW_\alpha(M)$ (here we use the assumption on $k$).
This defines an element of
$\tCO^a$ and completes the proof of the proposition.
\end{proof}

We denote by $p_a:\tCO\to\tCO^a$ the projection functor.

\subsubsection{}
\label{constructionproj}

We now give a construction of projective objects (under Hypothesis \ref{width}).

\begin{lemma}
Let $a\in\CP$ and $d\in\BZ$.
There is an integer $r$ such that the canonical map
$$\Hom(\tDelta_m(H)\langle -d\rangle,M)\to M_d$$ is an isomorphism
for all $m\ge r$ and $M\in\tCO^a$.
\end{lemma}

\begin{proof}
Replacing $M$ by $M\langle d\rangle$ and $a$ by $a+d$, we can assume
that $d=0$.

There is an integer $r$ such that
$p_a(\tDelta(H)\langle r'\rangle)=0$ for $r'\ge r$.
The exact sequence $$0\to \tDelta(B^m\otimes H)\langle m\rangle \to
\tDelta_m(H)\to \tDelta_{m-1}(H)\to 0$$
shows that the canonical map
$$\Hom(\tDelta_r(H),M)\iso\Hom(\tDelta_m(H),M)$$
is an isomorphism for any $M\in\tCO^a$ and $m\ge r$.
Equivalently, the canonical map
$$\Hom_B(B/B^{\ge r},M)\iso\Hom_B(B/B^{\ge m},M)$$
is an isomorphism. Since $M$ is locally $B$-nilpotent, this gives an
isomorphism $$\Hom(\tDelta_m(H),M)\iso M_0.$$
\end{proof}

\begin{cor}
\label{projectivegrad}
Let $E\in\Irr(H(k_0))$ and $a\in c_E+\BZ$. Then, the object
$p_a(\tDelta_r(k\otimes_{k_0} E)\langle a-c_E\rangle)$ of $\tCO^a$ is
independent of $r$,
for $r\gg 0$. It is
projective, has a filtration by modules
$\tDelta(k\otimes_{k_0} F)\langle r\rangle$ and has a quotient isomorphic to
$\tDelta(k\otimes_{k_0} E)\langle a-c_E\rangle$.
\end{cor}

\begin{cor}
\label{projective}
Let $E\in\Irr(H(k_0))$. Then, for $r\gg 0$, the module
$\Delta_r(k\otimes_{k_0} E)$
has a projective direct summand which is $\Delta$-filtered and
has a quotient isomorphic to $\Delta(k\otimes_{k_0} E)$.
\end{cor}

\begin{cor}
\label{progenerator}
There is an integer $r$ such that $\Delta_r(H)$ contains a progenerator
of $\CO$ as a direct summand.
\end{cor}

\begin{lemma}
\label{Ext1Delta}
Let $E,F\in\Irr(H(k_0))$ such that
$\Ext^1_\CO(\Delta(k\otimes_{k_0} E),\Delta(k\otimes_{k_0} F))\not=0$. Then,
$c_F-c_E$ is a positive integer.
\end{lemma}

\begin{proof}
By Lemma \ref{actionE} and Proposition \ref{decomposition},
we have $\Ext^1_\CO(\Delta(k\otimes_{k_0} E),\Delta(k\otimes_{k_0} F))=0$
if $c_F-c_E$ is not an integer. Assume now $c_F-c_E$ is an integer. Then
\begin{align*}
\Ext^1_\CO(\Delta(k\otimes_{k_0} E),\Delta(k\otimes_{k_0} F))&\simeq
\Ext^1_{\tCO}(\tDelta(k\otimes_{k_0} E),\tDelta(k\otimes_{k_0} F)
\langle c_F-c_E\rangle)\\
&\simeq
\Ext^1_A(\tDelta(k\otimes_{k_0} E),\tDelta(k\otimes_{k_0} F)
\langle c_F-c_E\rangle),
\end{align*}
by Lemma \ref{actionE} and Proposition \ref{decomposition}.
Now, 
$$\Ext^1_A(\tDelta(k\otimes_{k_0} E),\tDelta(k\otimes_{k_0} F)
\langle c_F-c_E\rangle) \simeq
\Ext^1_{BH}(k\otimes_{k_0} E,\Res_{BH}\tDelta(k\otimes_{k_0} F)
\langle c_F-c_E\rangle).$$
If the last $\Ext^1$ is non zero, then
$c_F-c_E$ is a positive integer.
\end{proof}

\begin{cor}
\label{goodproj}
Assume $k$ is a field.
Let $E\in\Irr(H)$. Then, $L(E)$ has a projective cover
$P(E)$ with a filtration 
$Q_0=0\subset Q_1\subset\cdots\subset Q_d=P(E)$ such that
$Q_i/Q_{i-1}\simeq \Delta(F_i)$ for some $F_i\in\Irr(H)$,
$c_{F_i}-c_E$ is a positive integer for $i\not=d$ and $F_d=E$.
\end{cor}

\begin{proof}
We know already that there is an  indecomposable projective 
module $P(E)$ as in the statement satisfying
all assumptions but the one on $c_{F_i}-c_E$, by
Corollary \ref{projective}.

Take $r\not=d$ maximal such that
$Q_r/Q_{r-1}\simeq \Delta(F)$ with $c_F-c_E$ not a positive integer.
By Lemma \ref{Ext1Delta}, the extension of $P(E)/Q_{r-1}$ by
$\Delta(F)$ splits. So, we have a surjective morphism
$P(E)\to \Delta(E)\oplus \Delta(F)$. This is impossible since $P(E)$ is 
indecomposable and projective.
\end{proof}

\pagebreak[3]\subsection{Highest weight theory}
\subsubsection{}
We assume here that $k$ is a field.

For $E$ a simple $H$-module,
all proper submodules of $\Delta(E)$ are graded submodules by
Proposition \ref{decomposition}, hence are contained in
$\Delta(E)_{>0}$. Consequently, $\Delta(E)$ has a unique maximal
proper submodule, hence a unique simple quotient which we denote by $L(E)$.

It follows from (\ref{Extvanishing}) that $L(E)$ is the unique
simple submodule of $\nabla(E)$ and that $L(E)\not\simeq L(F)$ for
$E\not\simeq F$.

\begin{prop}
\label{simpleobjects}
The simple objects of $\COln$ are the $L(E)$ for $E\in\Irr(H)$.
\end{prop}

\begin{proof}
Let $N\in\COln$. Then there is a simple
$H$-module $E$ such
that $\Hom_{BH}(E,\Res_{BH}N)\not=0$. By (\ref{adjunctionDelta}), it follows
that every simple object of $\COln$ is a quotient of
$\Delta(E)$ for some simple $H$-module $E$.
\end{proof}

\subsubsection{}
Let $M$ be a $A$-module.
Let $p(M)$ be the set of elements of $M$ annihilated by $B^{>0}$. This is
an $H$-submodule of $M$.

\begin{lemma}
\label{quotientDelta}
Let $M$ be a $A$-module and $E$ an $H$-module.
Then,
\begin{itemize}
\item
$M$ is a quotient of $\Delta(E)$ if and only if there is a
morphism of $H$-modules $\varphi:E\to p(M)$ such that
$M=A \varphi(E)$ ;
\item
If $k$ is a field and $E$ is simple, then $M\simeq L(E)$ if
and only if $M=Ap(M)$ and $p(M)\simeq E$.
\end{itemize}
In particular, $Ap(M)$ is the largest submodule of $M$ that is
a quotient of $\Delta(F)$ for some $H$-module $F$.
\end{lemma}

\begin{proof}
The first assertion follows from (\ref{adjunctionDelta})
and the isomorphism 
$$\Hom_{BH}(E,\Res_{BH}M)\simeq\Hom_H(E,p(M)).$$

Now, we assume $k$ is a field and $E$ is simple.

Assume $p(M)\simeq E$ and $M=Ap(M)$. Then, $M$ is in $\COln$.
Let $N$ be a non-zero submodule of $M$. We have
$0\not= p(N)\subseteq p(M)$, hence $p(N)=p(M)$ and $N=M$. So, $M$ is simple
and isomorphic to $L(E)$ since $M$ is a quotient of $\Delta(E)$.

Assume $M\simeq L(E)$. Since $\dim_k\Hom(\Delta(F),L(E))=1$ if
$E\simeq F$, and this Hom-space  vanishes otherwise,
it follows from
(\ref{adjunctionDelta}) that $p(M)\simeq E$.
\end{proof}

Let $M\{0\}=0$ and define by induction
$N\{i\}=M/M\{i\}$, $L\{i\}=Ap(N\{i\})$ and $M\{i+1\}$ as the inverse
image of $L\{i\}$ in $M$.
We have obtained a sequence of submodules of $M$,
$0=M\{0\}\subset M\{1\}\subset\cdots \subset M$.

Since $\Delta(E)$ is locally nilpotent for $B$,
the following proposition is clear. It describes how the objects of $\COln$
are constructed from $\Delta(E)$'s (cf Proposition \ref{charlocnilp}).

\begin{prop}
A $A$-module $M$ is locally nilpotent for $B$ if and only if
$\bigcup_i M\{i\}=M$, {\em i.e.}, if $M$ has a
filtration whose successive quotients are quotients of $\Delta(E)$'s.
\end{prop}

\begin{lemma}
\label{jordanfinite}
Assume $k$ is a field.
Every $A$-module quotient $M$ of $\Delta(E)$ has a finite Jordan-H\"older
series $0=M^0\subset M^1\subset\cdots\subset M^d=M$
with quotients $M^i/M^{i-1}\simeq L(F_i)$ such that $F_i\in\Irr(kW)$,
$c_E-c_{F_i}$ is a positive integer for $i\not=d$ and $F_d=E$.
\end{lemma}

\begin{proof}
By Proposition \ref{decomposition}, we can assume $c_F-c_E\in\BZ$.
  From Proposition \ref{decomposition}, it follows that $M$ inherits a grading
from $\Delta(E)$ (with $M_0\not=0$ and $M_{<0}=0$).
Note that since $\partial'$ acts as zero on $p(M)$, we have
$p(M)\subseteq\bigoplus_{F\in\Irr(H)}M_{c_E-c_F}$.

We will first show that $M$ has a simple submodule.

Take $i$ maximal such that $p(M)_i\not=0$ and $F$ a simple
$H$-submodule of $p(M)_i$. Let $L=AF$. Then,
$p(L)\subseteq p(M)\cap M_{\ge i}$ and $p(L)\subseteq F\oplus L_{>i}$, hence $p(L)=F$.
It follows from Lemma \ref{quotientDelta} that $L\simeq L(F)$ and we are
done.

\smallskip
Let $d(M)=\sum_{F\in\Irr(H)}\dim M_{c_E-c_F}$.

We put $M'=M/L$. We have $d(M')<d(M)$. So, the lemma follows
by induction on $d(M)$.
\end{proof}

\pagebreak[3]\subsection{Properties of category $\CO$}
\label{propcatO}
We assume here in \S \ref{propcatO} that $k$ is a field.
We now derive structural properties of our categories.
\subsubsection{}

\begin{cor}
\label{lfJordan}
Every object of $\COln$ has an ascending filtration 
whose successive quotients are semi-simple.
\end{cor}

\begin{proof}
Follows from Lemma \ref{jordanfinite} and Proposition \ref{charlocnilp}.
\end{proof}

\begin{cor}
\label{artinian}
Every object of $\CO$ has a finite Jordan-H\"older series.
\end{cor}

\begin{proof}
The multiplicity of $L(E)$ in a filtration of $M\in\CO$ given by
Corollary \ref{lfJordan} is bounded by $\dim \CW_{-c_E}(M)$, hence
the filtration must be finite.
\end{proof}

\begin{cor}
The category $\tCO^a$ is generated by the
$L(E)\langle r\rangle$, with $r=c_E-a$.
\end{cor}

\begin{proof}
Follows from Lemma \ref{actionE} and Proposition \ref{simpleobjects}.
\end{proof}

\begin{cor}
\label{decomposition2}
Given $a\in k$, the full abelian Serre subcategory of the category of
$A$-modules generated by the $L(E)$ with $c_E\in a+\BZ$ is
$\CO^{a+\BZ}$.
\end{cor}

\subsubsection{}\label{reciproc}
\begin{thm}
\label{highestweightcategory}
The category $\CO$ is a highest weight category (in the sense of
\cite{CPS1})
with respect to the
relation: $E< F$ if $c_F-c_E$ is a positive integer.
\end{thm}

\begin{proof}
Follows from Corollary \ref{goodproj} and Lemma \ref{jordanfinite}.
\end{proof}

The standard and costandard objects are the $\Delta(E)$ and $\nabla(E)$.
There are projective modules $P(E)$, injective modules $I(E)$,
tilting modules $T(E)$. We have reciprocity formulas, cf.
\cite[Theorem 3.11]{CPS1}:
$$[I(E):\nabla(F)]=[\Delta(F):L(E)]\quad\mbox{and}
\quad [P(E):\Delta(F)]=[\nabla(F):L(E)].
$$

\begin{cor}
\label{semisimple}
If $c_E-c_F\not\in\BZ-\{0\}$ for all $E,F\in\Irr(H(k_0))$, then
$\CO$ is semi-simple.
\end{cor}

\begin{prop}
Let $M\in\CO$.
The following assertions are equivalent
\begin{itemize}
\item
$M$ has a $\Delta$-filtration
\item
$\Ext^i_\CO(M,\nabla(H))=0$ for $i>0$
\item
$\Ext^1_\CO(M,\nabla(H))=0$
\item
the restriction of $M$ to ${\Bbar}$ is free.
\end{itemize}
\end{prop}

\begin{proof}
The equivalence between the first three assertions is classical. The
remaining equivalences follow from the isomorphism
$\Ext^1_{\CO}(M,\nabla(H))\iso\Ext^1_A(M,\nabla(H))\iso\Ext^1_{\Bbar}(M,k)$.
\end{proof}

\subsubsection{}
 From Proposition \ref{charlocnilp}, we deduce

\begin{lemma}
\label{finiteness}
Let $M\in\COln$. The following conditions are equivalent
\begin{itemize}
\item $M\in\CO$
\item $M$ is finitely generated as a $A$-module
\item $M$ is finitely generated as a $\Bbar$-module.
\end{itemize}
\end{lemma}

\begin{lemma}
There is $r\ge 0$ such that for all $M\in\CO$, $a\in k$ and
$m$ in the generalized eigenspace for $\partial'$ for the eigenvalue $a$,
then $(\partial'-a)^r m=0$.
\end{lemma}

\begin{proof}
The action of $\partial'$ on $\Delta(H)$ is semi-simple. It follows
that, given $r\ge 0$, $a\in k$ and $m\in\Delta_r(H)$ in the 
generalized eigenspace for $\partial'$ for the eigenvalue $a$, then
 $(\partial'-a)^{r} m=0$.

Now, by Corollary \ref{projective}, there is some integer $r$ such that
every object of $\CO$ is a quotient of $\Delta_r(H)^l$ for some $l$.
\end{proof}

\begin{prop}
There is $r\ge 0$ such that
every module in $\COln$ is generated by the kernel of~$B^{\ge r}$.
Further,
there is an integer $r>0$ such that for $M\in\COln$, we have
$M\{i\}=M\{r\}$ for $i\ge r$.
\end{prop}

\begin{proof}
Let $r\ge 0$ such that every projective indecomposable object in $\CO$
is a quotient of $\Delta_{r-1}(H)$. This means that every object in
$\CO$ is generated by the kernel of $B^{\ge r}$. Now, consider $M\in\COln$
and $m\in M$. Let $N$ be the $A$-submodule of $M$ generated by $m$.
This is in $\CO$, hence $m$ is in the submodule of $N$ generated 
by the kernel of $B^{\ge r}$.
\end{proof}

\begin{prop}
Every object in $\COln$ is generated by the $0$-generalized
eigenspace of $\partial'$.
\end{prop}

\begin{proof}
It is enough to prove the proposition for projective indecomposable
objects in $\CO$, hence for $\Delta_r$'s, where it is obvious.
\end{proof}
\subsubsection{}
\label{algebraforO}
Let $Q$ be a progenerator for $\CO$ (cf Corollary \ref{progenerator})
and $\Gamma=(\End_A Q)^{\opp}$. Then, $\Gamma$ is a finitely generated
projective $\CO$-module. 
We have mutually inverse standard equivalences
\begin{equation}\label{progen}\Hom(Q,-):\CO\iso \Gamma\mMod,
\ \ \ Q\otimes_\Gamma(-):\Gamma\mMod\iso\CO.\end{equation}

Let now $X$ be a (non-necessarily finitely generated) $\Gamma$-module.
Then, $Q\otimes_\Gamma X$ is a quotient of $Q^{(I)}$ for some set $I$, where
$X$ is a quotient of $\Gamma^{(I)}$. Now, $Q^{(I)}$ is in $\COln$. So,
the functor $Q\otimes_\Gamma (-):\Gamma\mMOD\to A\mMOD$ takes values in
$\COln$ and we have equivalences
$$\Hom(Q,-):\COln\iso \Gamma\mMOD,\ \ \ Q\otimes_\Gamma-:\Gamma\mMOD\iso\COln.$$

\section{Rational Cherednik algebras}
\subsection{Basic definitions}
\label{sectionreminder}
Let $V$ be a finite dimensional vector space and
$W\subset \GL(V)$ a finite complex reflection group.
Let 
$\CA$ be the set of reflecting hyperplanes of $W$. 
Given $H\in \CA,$  let $W_H\subset W$ be the subgroup formed by the
elements of $W$ that fix $H$ pointwise.
We choose $v_H\in V$ such that $\BC v_H$ is a $W_H$-stable
complement to $H$. Also, let $\alpha_H\in V^*$ be a linear form with
kernel $H$.

Let $k$ be a noetherian commutative $\BC$-algebra.
The group
$W$ acts naturally  on $\CA$ and on the group algebra $kW$,
by conjugation. Let $\gamma: \CA\to kW\,,\,H\mapsto \gamma_H,$
be a $W$-equivariant map such that
 $\gamma_H$ is an element of 
$kW_H\subset kW$ with trace zero, for each $H\in \CA$.

Given $\gamma$ as above,
one  introduces an associative $k$-algebra  $\sa(V,\gamma),$
the rational Cherednik algebra. It is defined
as  the quotient of  $k\otimes_\BC T(V\oplus V^*)\rtimes W$,
the cross-product of $W$ with $k$-tensor algebra,
by the relations
$$[\xi,\eta]=0 \text{ for } \xi,\eta\in V,\ \ 
[x,y]=0 \text{ for } x,y\in V^*$$
$$[\xi,x]=\langle \xi,x\rangle+\sum_{H\in\CA}
\frac{\langle\xi,\alpha_H\rangle\langle v_H,x\rangle}
     {\langle v_H,\alpha_H\rangle} \gamma_H$$

\begin{rem} Let ${\operatorname{Refl}}\subset W$  denote the set of (pseudo)-reflections.
Clearly ${\operatorname{Refl}}$ is an ${\operatorname{Ad}} W$-stable subset.
Giving $\gamma$  as above is equivalent to giving 
 a $W$-invariant function $c: {\operatorname{Refl}} \to k\,,\, g\mapsto c_g$
such that $\gamma_H= \sum_{g\in W_H\smallsetminus\{1\}}\,c_g\cdot g$.
One may use the function $c$ instead of $\gamma$, and 
write   $v_g\in V$, resp. $\alpha_g\in V^*$, instead of $v_H$, resp.
$\alpha_H$,
for any
$g\in W_H\smallsetminus\{1\}$. Then
 the last commutation relation in the algebra $\sa(V,\gamma)$
reads:
$$[\xi,x]=\langle \xi,x\rangle+\sum_{g\in {\operatorname{Refl}}}
c_g \cdot\frac{\langle\xi,\alpha_g\rangle\langle v_g,x\rangle}
     {\langle v_g,\alpha_g\rangle}\cdot g,
$$
which is essentially the commutation relation used in \cite{EtGi}.
In case of a Weyl group $W$, in \cite{EtGi,BeEtGi,Gu}, the
coefficients $c_\alpha$ ($\alpha$ a root) were used instead of the
$c_g$'s.
Then, $\gamma_H=-2c_\alpha g$ for $H$ the kernel
of $\alpha$ and $g$ the associated reflection.
\end{rem}
\begin{rem}\label{eps} Put $e_H=|W_H|$.
Denote by $\eps_{H,j}=\frac{1}{e_H}\sum_{w\in W_H}\det(w)^jw$ the
idempotent of $\BC W_H$ associated to the character $\det^{-j}_{|W_H}$.
Given $\gamma$ as above, there is a unique family
$\{k_{H,i}=k_{H,i}(\gamma)\}_{H\in\CA/W\,,\,0\le i\le e_H}$ of elements of $k$ 
such that
 $k_{H,0}=k_{H,e_H}=0$ and
$$\gamma_H=e_H\sum_{j=0}^{e_H-1}
(k_{H,j+1}(\gamma)-k_{H,j}(\gamma))\eps_{H,j}.$$
We observe that $\gamma$ can be recovered from the $k_{H,i}(\gamma)$'s
by the formula
$$\gamma_H=\sum_{w\in W_H-\{1\}}
\left(\sum_{j=0}^{e_H-1}\det(w)^j\cdot(k_{H,j+1}(\gamma)-k_{H,j}(\gamma))\right)w.$$
This way, we get back to the definition of \cite{DuOp}. 
\end{rem}

Introduce free
commutative positively graded $k$-algebras $\sp=k\otimes_\BC S(V^*)=\bigoplus_{i\geq 0}\sp_i$ and
${\ss}=k\otimes_\BC S(V)=\bigoplus_{i\geq 0}\ss^i$.
We have a triangular decomposition
$\sa=\sp\otimes_k  kW\otimes_k {\ss}$ as $k$-modules
\cite[Theorem 1.3]{EtGi}.

For $H\in\CA$, we put
$$a_H(\gamma)=\sum_{i=1}^{e_H-1} e_H\cdot k_{H,i}(\gamma)\cdot
\eps_{H,i}\in k[W_H]
\quad\text{and}\quad
z(\gamma)=\sum_{H\in\CA} a_H(\gamma)\in Z(kW).$$

We denote by $\Irr(kW)$ a complete set of representatives of
isomorphism classes of simple $kW$-modules.

For $E\in\Irr(\BC W)$, we denote by $c_E=c_E(\gamma)$ the scalar by which $z(\gamma)$ acts
on $k\otimes_\BC E$.
The elements $a_H(\gamma)\,,\,z(\gamma),$ and $c_E(\gamma),$
may be thought of
as functions of the coefficients $k_{H,i}=k_{H,i}(\gamma)$
(through their dependence on $\gamma$).
In particular, it was shown in \cite[Lemma 2.5]{DuOp}
that $c_E$, expressed as a function of the   $k_{H,i}$'s,
is a linear function with non-negative integer
coefficients.
\smallskip

Below, we will often use simplified notations and write
$\sa$ for $\sa(\gamma)$, $\,k_{H,i}$ for $k_{H,i}(\gamma)$, and
$\,z$ for $z(\gamma),$ etc.

We introduce a grading on $\sa$ by putting $V^*$ in degree $1$, $V$ in degree
$-1$ and $W$ in degree $0$ (thus, the induced grading on the subalgebra
$\sp\subset \sa$ coincides with the standard one on $\sp$, while
the  induced grading on the subalgebra
$\ss\subset \sa$  {\it  differs by a sign} from the standard 
one on~$\ss$).\q{We use superscripts to indicate
 the standard (non-negative) grading on $\ss$, and subscripts to
denote the gradings on $\sa$ and $\sp$. Thus, putting formally
$\ss_{-i}:=\ss^i$ one recovers compatibility: $\ss_i=\ss\cap\sa_i$.}

Let $\se_k=\sum_{b\in\CB}b^\vee b$ be the ``deformed Euler vector field'',
where $\CB$ is a basis of $V$ and $\{b^\vee\}_{b\in\CB}$ is the dual
basis. We also put $\se=\se_k-z$. 
The elements $\se_k$ and $\se$ commute with $W$. Note that
$\sum_b\, [b,b^\vee]=\dim V+\sum_H \gamma_H$.

We have
\begin{equation}
\label{relationspartial}
[\se,\xi]=-\xi \text{ and } [\se,x]=x \text{ for }
\xi\in V \text{ and } x\in V^*.
\end{equation}

This shows the grading on $\sa$ is ``inner'', \ie,
$\sa_i=\{a\in \sa\; |\; [\se,a]=i\cdot a\}$.

\pagebreak[3]
\subsection{Category $\CO$ for the rational Cherednik algebra}
\label{catOCherednik}
\subsubsection{}
We apply now the  results of \S\ref{sectioncatO} in the special case:
$A=\sa=\sa(V,\gamma)$, $B=\ss$, $\Bbar=\sp$, $H=kW$, $k_0=\BC$, $H(k_0)=\BC W$, 
$\partial=\se$, $\partial'=\se_k$ and $\partial_0=z$.
In particular, we have the category $\CO(\gamma):=
\CO(\sa(V,\gamma))$, which was 
first considered, in the setup of Cherednik algebras, in \cite{DuOp}.

For any (commutative) algebra map $\psi:k\to k'$,
there is  a base extension
functor $\CO(\gamma)\to \CO(\psi(\gamma))$ given by
$\sa(\psi(\gamma))\otimes_\sa (-)$.

\subsubsection{}
Assume $k$ is a field. Since $\CO$ and $\tCO$ have finite global
dimension (Theorem \ref{highestweightcategory}),
the Grothendieck group of the category of modules coincides
with the Grothendieck group $K_0$ of projective modules.

We have a morphism of $\BZ$-modules
$f:K_0(\tCO)\to \BZ[[q]][q^{-1}]\otimes K_0(\BC W)$
given by taking the graded character of the restriction of the module to $W$~:
$$M\mapsto \sum_{E\in\Irr(kW)}\sum_i q^i \dim\Hom_{kW}(E,M_i)[E].$$
Set $[\sp]:= \sum_{E\in\Irr(kW)}\sum_i q^i
\dim\Hom_{kW}(E,{\sp}_i)\cdot[E]$.
This is an invertible element of $\BZ[[q]][q^{-1}]\otimes K_0(\BC W)$,
and for any $F\in\Irr(kW)$, we
have $f([\Delta(F)])=[\sp]\cdot [F].$ Since the classes of standard modules
generate the $K_0$-group, 
we obtain an isomorphism
$\frac{1}{[\sp]} f:K_0(\tCO)\iso \BZ[q,q^{-1}]\otimes K_0(\BC W)$.
 
\smallskip
Let $k[(\Bk_{H,i})_{1\le i\le e_H-1}]$ be the polynomial ring in the
indeterminates $\Bk_{H,i}$ with $\Bk_{w(H),i}=\Bk_{H,i}$ for $w\in W$.
We have a canonical evaluation morphism $k[(\Bk_{H,i})]\to k$ given by the choice of
parameters.
Let $\Gm$ be the kernel of that morphism, $R$ the completion of
$k[(\Bk_{H,i})]$ at  $\Gm$, and $K$ the field of fractions of
$R$.

We have a decomposition map $K_0(\CO_K)\iso K_0(\CO)$. It
sends $[\Delta(K\otimes_\BC E)]$ to $[\Delta(k\otimes_\BC E)]$.

\begin{prop}
\label{classdeltanabla}
Assume $k$ is a field. Then, 
$[\Delta(E)]=[\nabla(E)]$ and $[P(E)]=[I(E)]$ for any $kW$-module $E$.
\end{prop}

\begin{proof}
We first consider the equality $[\Delta(E)]=[\nabla(E)]$.
The corresponding statement for $K$ is true, since the category is
semi-simple in that case. Hence the modules are determined, up to
isomorphism, by their socle (resp. by their head).

The statement for $k$ follows by using the decomposition map.

Now, the equality $[P(E)]=[I(E)]$ follows, using the reciprocity formulas
(\S \ref{reciproc}).
\end{proof}

\pagebreak[3]
\section{Duality, Tiltings, and Projectives}
\label{sectionDmod}
\subsection{Ringel duality}
We keep the setup of \S\ref{triag}, with $k$ being a field.
We
 make the following two additional assumptions
\begin{itemize}\label{cond2}
\item We have $\Bbar\otimes H \otimes B= B\otimes H \otimes\Bbar = A$;
\item The
subalgebra $B\subset A$ is {\em Gorenstein} (with parameter $n$), i.e.,
there exists an integer $n$  such that
$$
\Ext_B^i(k,B)=\begin{cases}k&\text{if }\; i=n\\
0&\text{if }\; i\neq n.
\end{cases}
$$
\end{itemize}
The Gorenstein condition implies that, for any
$E\in \Irr(H)$, viewed as a $BH$-module via the projection
$BH\to H$, we have $\Ext_{BH}^i(E, BH)= 0,$ for all
$i\neq n$; moreover,
$\Ext_{BH}^n(E, BH)= E^\flat,$
where $E^\flat$ is a {\em right} $BH$-module such that
$\dim E^\flat=\dim E$.

Assume further
 that the algebra
$A$ has finite homological dimension.
Thus (see \cite{Bj}), there is a well-defined
 duality functor
$$R\Hom_A(-,A):D^b(A\mMod)\iso D^b(A^\opp\mMod)^\opp.$$
Furthermore, this functor is an equivalence with inverse
$R\Hom_{A^\opp}(-,A)$.

The triangular decomposition $A=B\otimes H \otimes\Bbar$  gives a similar decomposition
$A^\opp={\Bbar}^\opp\otimes H^\opp\otimes B^\opp$,
for the opposite algebras.
Therefore, we may consider the category $\CO(A^\opp)$
and, for any simple  {\em right} $H$-module $E'$,
introduce the standard $A^\opp$-module
$$\Delta^\opp(E'):=\Ind_{(BH)^\opp}^{A^\opp}E'
=E'\otimes_{BH} A,
$$
 and also
the projective object $P^\opp(E')\in\CO(A^\opp)$, the tilting object
$T^\opp(E')\in\CO(A^\opp)$, etc.

\begin{lemma}\label{key}
The functor $R\Hom_A(-,A[n])$
sends $\Delta(E)$ to $\Delta^\opp(E^\flat),$ for
$E$ a finite-dimensional $H$-module.
\end{lemma}

\begin{proof}
Using that $A$ is free as a left $BH$-module, we compute
$$\Ext^i_A(\Delta(E),A)\iso
\Ext^i_{{B}H}(E,A)\iso \Ext^i_{{B}H}(E,BH)\otimes_{{B}H} A.$$
We see that this space vanishes for $i\neq n$, and for $i=n$ we get
$R\Hom_A\bigl(\Delta(E)\,,\,A[n]\bigr)\simeq
E^\flat\otimes_{{B}H} A=\Delta^\opp(E^\flat)$.
\end{proof}

We would like to use the duality functor $R\Hom_A(-,A[n])$  to obtain a functor
$D^b(\CO(A))\to D^b(\CO(A^\opp))^\opp$.
To this end, we will exploit a general result below
 valid for arbitrary highest weight categories
(a contravariant version of Ringel duality \cite[\S 6]{Ri}).

Given an additive category $\mathcal{C}$, let $K^b(\mathcal{C})$ be
the
corresponding homotopy category of bounded complexes in
$\mathcal{C}$.

\begin{prop}
\label{Ringel}
Let $A$ and $A'$ be two quasi-hereditary algebras and
$\CC=A\mMod$, $\CC'=A'\mMod$ the associated highest weight categories.
Let $F$ be a contravariant equivalence between the exact categories of 
$\Delta$-filtered objects
$F:\CC^\Delta\iso (\CC^{\prime\Delta})^\opp$.
Then,

\noindent
$\bullet\quad$
$F$ restricts to equivalences
$$\CC\mproj\iso (\CC'\mtilt)^\opp\quad\text{ and }\quad
\CC\mtilt\iso (\CC'\mproj)^\opp$$

\noindent
$\bullet\quad$ The canonical equivalences
$K^b(\CC\mproj)\iso D^b(\CC)$ and
$K^b(\CC'\mtilt)\iso D^b(\CC')$,
yield an equivalence of derived categories
$$D:D^b(\CC)\iso D^b(\CC')^\opp.$$

\noindent
$\bullet\quad$ 
Let $T=F(A)$, an $(A\otimes A')$-module. Then, 
we have $D=R\Hom_A(-,T)$ and $D^{-1}=R\Hom_{A'}(-,T)$.
Via duality $(A\mMod)\iso (A^\opp)\mMod^\opp$, the functor $D$ identifies
$\CC^\opp$ with the Ringel dual of $\CC'$.
\end{prop}

\begin{proof}
Let $M\in\CC^\Delta$. Then, 
$M$ is projective if and only if $\Ext^1(M,\Delta(E))=0$
for every standard object $\Delta(E)$ of $\CC$ (indeed, if
$0\to M'\to P\to M\to 0$ is an exact sequence with $P$ projective, then
$M'$ is $\Delta$-filtered, hence the sequence splits).
The module $F(M)$ is tilting if and only if 
$\Ext^1(\Delta(E'),F(M))=0$ for every  standard object $\Delta(E')$ of $\CC'$.
We deduce that $M$ is projective if and only if $F(M)$ is tilting.

So, $F$ restricts to equivalences $\CC\mproj\iso (\CC'\mtilt)^\opp$ and
$\CC\mtilt\iso (\CC'\mproj)^\opp$.

The last assertions of the Proposition are clear.
\end{proof}

We can now apply this construction to the category $\CO(A)$.
Specifically, Lemma \ref{key} implies that the
 functor $R\Hom_A(-,A[n])$ restricts to an equivalence
$\CO(A)^\Delta\iso (\CO(A^\opp)^\Delta)^\opp$.
Therefore, using Proposition \ref{Ringel} we immediately
obtain the following
\begin{prop}
\label{dualityequivalences}
The functor $R\Hom_A(-,A[n])^\opp$
restricts to equivalences
$$\CO(A)\mproj\iso (\CO(A^\opp)\mtilt)^\opp
\quad\text{ and }\quad
\CO(A)\mtilt\iso (\CO(A^\opp)\mproj)^\opp\,.$$
%\medskip

The canonical equivalences
$K^b(\CO(A)\mproj)\iso D^b(\CO(A))$ and
$K^b(\CO(A^\opp)\mtilt)\iso D^b(\CO(A^\opp))$
induce an equivalence
$D:D^b(\CO(A))\iso D^b(\CO(A^\opp))^\opp,$
such that $\Delta(E)\mapsto \Delta^\opp(E^\flat)\,,\,P(E)\mapsto
T^\opp(E^\flat),$ and $T(E)\mapsto P^\opp(E^\flat).$\qed
\end{prop}
\begin{cor}\label{ringel_cor}
The category $\CO(A^\opp)^\opp$ is the Ringel dual of $\CO(A)$.\qed
\end{cor}

\subsection{Naive duality for Cherednik algebras}\label{duality1}
 Recall the setup
of \S\ref{catOCherednik}.

Denote by $(-)\ii:\BC W\iso \BC W$ the anti-involution given by
$w\mapsto w\ii:=w^{-1}$ for $w\in W$. 

In this section, we compare the algebras $\sa
=\sa(V,\gamma)$
 and $\sa(V^*,\gamma\ii)$. This will provide us with means 
to switch between left and right modules, between
${\ss}$-locally finite and ${\sp}$-locally finite modules.

\smallskip
The anti-involution  $(-)\ii:\BC W\iso \BC W$ extends to an isomorphism
\begin{equation}\label{dag}
(-)\ii: \sa(\gamma)\iso \sa(\gamma\ii)^\opp\quad,\quad
V\ni\xi\mapsto -\xi,\ \ V^*\ni x\mapsto x,\ \  W\ni w\mapsto w^{-1}.
\end{equation}

\begin{rem}
If all pseudo-reflections of $W$ have order $2$, then
$\gamma\ii=\gamma$.
\end{rem}

Further, we define an isomorphism of $k$-algebras reversing the gradings by
$$\varphi:\sa(V,\gamma)\iso \sa(V^*,\gamma\ii)^{\opp}$$
$$V\ni\xi\mapsto \xi, \ \
V^*\ni x\mapsto x,\ \ 
W\ni w\mapsto w^{-1}.$$

\begin{rem}
When $V$ is self-dual, an isomorphism of $\BC W$-modules
$F:V\iso V^*$ extends to an algebra isomorphism (Fourier transform)
$$F:\sa(V,\gamma)\iso \sa(V^*,\gamma)$$
$$V\ni\xi\mapsto F(\xi), \ \
V^*\ni x\mapsto -F^{-1}(x),\ \ 
W\ni w\mapsto w.$$
The functor $F_*$ restricts to an equivalence $\CO(V,\gamma)\iso \CO(V^*,\gamma\ii)$.
\end{rem}

\subsubsection{}
Given $M\in \COln$, denote by $M^\vee$ the $k$-submodule of
${\sp}$-locally nilpotent elements of $\Hom_k(M,k)$. This is a right
$\sa$-module. Via $\varphi_*$, this becomes a left
$\sa(V^*,\gamma\ii)$-module.
If $M$ is graded, then $M^\vee=\Homgr_k^\bullet(M,k)$.

Thus we have defined a functor (analogous of the standard duality on
the category $\CO$ in the Lie algebra case):
\begin{equation}\label{vee}
(-)^\vee:\COln(V,\gamma)\to \COln(V^*,\gamma\ii)^\opp
\end{equation}
When $k$ is a field, this functor is an equivalence.

Given a $kW$-module $E$, we use the notation $E^\vee=\Hom_k(E,k)$ for the
dual $kW$-module.

\begin{prop}
\label{usualduality}
We have $\Delta(E)^\vee\iso \nabla(E^\vee)$ for any $kW$-module $E$.
If $k$ is a field, then
$$ L(E)^\vee\iso L(E^\vee)\;,\; P(E)^\vee\iso
I(E^\vee)\;,
\; I(E)^\vee\iso P(E^\vee)\;,
\; 
\nabla(E)^\vee\iso \Delta(E^\vee)
\;,\;
T(E)^\vee\iso T(E^\vee).
$$
\end{prop}

\begin{proof}
We have
$$\Homgr^\bullet_k(\sa\otimes_{{\ss}W}E,k)\iso
\Homgr^\bullet_{({\ss}W)^\opp}(\sa,\Hom_k(E,k))$$
and the first part of the proposition follows.

The second assertion follows from the characterization of $L(E)$
(resp. $L(E^\vee)$) as the unique simple quotient (resp. submodule)
of $\Delta(E)$ (resp. $\nabla(E^\vee)$).
The other assertions are immediate consequences of the homological
characterizations of the objects and/or the existence of suitable filtrations.
\end{proof}

Note that the functor $(-)^\vee$ restricts to a functor $\CO(V,\gamma)\to \CO(V^*,\gamma\ii)^\opp$.
When $k$ is a field, it is an equivalence.
A compatible
choice of progenerators for $\CO(V,\gamma)$ and $\CO(V^*,\gamma\ii)$ gives then
an isomorphism between the algebra $\Gamma(V)$ for $\CO(V,\gamma)$
and the oppposite algebra $\Gamma(V^*)^\opp$ for $\CO(V^*,\gamma\ii)$
(cf \S \ref{algebraforO}).

\begin{cor}\label{eric_cor}
Let $E$ and $F$ be two simple $kW$-modules. Then,
the multiplicity of $\Delta(E)$ in a $\Delta$-filtration of
$P(F)$, for $\CO(V,\gamma)$, is equal to the multiplicity of
$L(F^\vee)$ in a composition series of
$\Delta(E^\vee)$, for $\CO(V^*,\gamma\ii)$.
\end{cor}

\begin{proof}
By \S \ref{highestweightcategory}, the
multiplicity of $\,\nabla(E^\vee)\,$ in a $\nabla$-filtration of
$\,I(F^\vee)\,$ is equal to the multiplicity of
$\,L(F^\vee)\,$ in a composition series of
$\,\Delta(E^\vee)$. 

The functor $(-)^\vee$ sends $P(F)$ to $I(F^\vee)$ and
$\Delta(E)$ to $\nabla(E^\vee)$ (Proposition \ref{usualduality})
and the result follows.
\end{proof}

\begin{rem}
When $k$ is a field and $W$ is real, we obtain, via Fourier transform,
a duality on $\COln$ and on $\CO$.
Since all complex representations of $W$ are
self-dual, we have then $\Delta(E)^\vee\iso \nabla(E)$.
\end{rem}

\pagebreak[3]
\subsection{Homological properties of Cherednik algebras}
\subsubsection{}
The rational Cherednik algebra is a deformation of the 
cross-product of $W$ with the Weyl algebra 
of polynomial differential operators on $V$. In particular, 
there is a standard increasing filtration on $\sa$ with $W$ placed in 
degree $0$ and $V\oplus V^*$
in degree $1$.  The associated graded ring,  $\gr\sa$, is isomorphic to
$S(V\oplus V^*) \rtimes W$ \cite[\S 1]{EtGi}.
It follows (see \cite{Bj}, \cite[\S V.2.2]{Bo}), that $\sa$
is left and right noetherian, provided $k$ is noetherian.
Since $V^*\oplus V$ is a smooth variety of
dimension $2\dim V$, the algebra $\sa$ has homological dimension at most $2\dim V$.
Furthermore, the usual results and concepts on $\CD$-modules
(characteristic variety, duality) also make sense for $\sa$,
even though the algebra $\gr\sa$ is not commutative.

\subsubsection{}

We assume $k$ is a field, and put $n=\dim V$.
The algebras $\ss$ and $\sp$ are
clearly Gorenstein with parameter $n$.
Moreover, we have $\Ext^n_{\ss}(k,{\ss})\simeq\Lambda^n V^*$.
Hence, $E^\flat=\Lambda^n V^*\otimes \Hom_k(E,k)=\Lambda^n V^*\otimes E^\vee,$
for any finite dimensional $W$-module $E$.

It will be useful to compose the functor
 $R\Hom_{\sa(\gamma)}(-,\sa(\gamma))$ with the anti-involution $(-)\ii$,
see \eqref{dag}, to get 
the following  composite equivalence
\begin{equation}\label{rdag}
R\Hom_{\sa(\gamma)}(-,\sa(\gamma))\ii:\,D^b(\sa(\gamma)\mMod)\,\iso\,
D^b(\sa(\gamma)^\opp\mMod)^\opp\,\iso\,
D^b(\sa(\gamma\ii)\mMod)^\opp.
\end{equation}

From Proposition \ref{Ringel} we immediately
obtain the following
\begin{prop}
\label{Cher_dualityequivalences}
The functor $R\Hom_\sa(-,\sa[n])\ii$ gives rise to an
 equivalence
$$D:D^b(\CO(\gamma))\iso D^b(\CO(\gamma\ii))^\opp.\qquad\qed$$
\end{prop}
%\smallskip

We further introduce an equivalence
$$(-)^\vee\circ D:\;D^b(\CO(V,\gamma))\iso
D^b(\CO(V^*,\gamma))$$
such that
$$\Delta(E)\mapsto \nabla(\Lambda^n V\otimes_\BC E)$$
$$T(E)\mapsto I(\Lambda^n V\otimes_\BC E)$$
$$P(E)\mapsto T(\Lambda^n V\otimes_\BC E)$$

In particular, we obtain (cf. Corollary \ref{ringel_cor})

\begin{cor}
The category $\CO(V^*,\gamma)$ is the Ringel dual of
$\CO(V,\gamma)$.\qed 
\end{cor}

\begin{rem}
Note that if $W$ is real, then $\CO$ is its own Ringel dual.
\end{rem}
\subsubsection{Semiregular bimodule.} Write
$\sp^{\circledast}=\oplus_i\,\Hom(\sp_i,k)$
for the graded dual of $\sp$, and form
the vector space $R:=\sp^{\circledast}\otimes_{k} \ss W$. Let
us fix an isomorphism of $\BC$-vector spaces
$\Lambda^n V\iso \BC$. We have the
following
canonical isomorphisms: 
$$ \Homgr^\bullet_{\ss W}(\sa,\Lambda^nV\otimes_{\BC}\ss W)
\stackrel{\sim}{\rightarrow}
\Homgr^\bullet_{k}(\sp,\ss W)\stackrel{\sim}{\leftarrow}
 \sp^{\circledast}\otimes_{k}\ss W\stackrel{\sim}{\rightarrow}\sp^{\circledast}\otimes_{\sp} \sa.$$
The first two isomorphisms
define a left $\sa$-module structure on $R$, and the last one
defines a right 
$\sa$-module structure on $R$. It is possible to check by explicit
calculations 
that the left and right $\sa$-module structures commute, 
so that $R$ becomes an $\sa$-bimodule.
It is a Cherednik algebra analogue of the {\em semiregular} bimodule,
considered in \cite{A}, \cite{So2} in the Lie algebra case.

From the isomorphisms  of $\sa$-modules
$$ R\otimes_{\sa} \Delta(E)=R\otimes_{\ss W} 
E\stackrel{\sim}{\rightarrow} \Homgr^\bullet_{\ss W}(\sa,\Lambda^nV\otimes E)
$$ 
we deduce
\begin{prop}
The functor $M\mapsto \Homgr^\bullet_{k}(R\otimes_{\sa} M,k)^\dag$ 
(= left $\sa(V,\gamma^{\dagger})$-module) sends
$\Delta(E)$ to $\Delta(\Lambda^nV^*\otimes_{\BC} E^{\vee})$. \qed
\end{prop}

\subsubsection{}
Given $M$ a finitely generated $\sa$-module, a good filtration of $M$
is a structure of filtered $\sa$-module on $M$ such that
$\gr M$ is a finitely generated $\gr \sa$-module. The characteristic
variety $\Ch(M)$ is the support of $\gr M$, viewed as a $W$-equivariant
sheaf
on $V^*\oplus V$ (a closed subvariety). It is well defined, i.e.,
 is independent of the choice of the good filtration (every
finitely generated $\sa$-module admits a good filtration).
Note that Bernstein's inequality: $\dim\Ch(M)\ge \dim V$ does not hold
in general. Further,  for $M$ in $\CO$, the complex
$D(M)$ has zero homology outside the degrees $0,\ldots,n$.

Let $T=\bigoplus_E T(E)$ where $E$ runs over the simple $kW$-modules.

\begin{cor}
Let $M\in\CO$. Then,
$\dim\Ch(M)=\dim V - \min\{i\;|\; \Ext^i_\CO(T,M)\neq 0\}$.
\end{cor}

\begin{proof}
Let $R=\End(T)^\opp$. The functor
$R\Hom_\CO(T,-):D^b(\CO)\iso D^b(R\mMod)$ is an equivalence.
Composing with the inverse of $(-)^\vee\circ D$ we obtain an
equivalence $D^b(\CO(V^*,\gamma))\iso D^b(R\mMod)$ that restricts
to an equivalence $\CO(V^*,\gamma)\iso R\mMod$.
We see that $\;\;\min\{i\;|\;\;$
$\Ext^i_\CO(T,M)\not=0\}$ $=\min\{i\;|\;
H^i(DM)\not=0\},$
where the RHS is equal to $\min\{i\;|\; \Ext^i_\sa(M,\sa)\not=0\}-\dim V$
by the definition of $D$.
The result now follows from the well-known formula,
see e.g \cite{Bj}:
$\dim\Ch(M)=2\dim V - \min\{i\;|\; \Ext^i_\sa(M,\sa)\not=0\}.$
\end{proof}

\pagebreak[3]
\section{Hecke algebras via monodromy}
\label{sectionHecke}

\pagebreak[3]\subsection{Localisation}
\subsubsection{}
\label{fullyfaithful}
Let $V\reg=V-\bigcup_{H\in\CA}H$ and ${\sp}\reg=k[V\reg]=
{\sp}[(\alpha_H^{-1})]_{H\in\CA}$. The algebra structure on
$\sa$ extends to an algebra structure on
$\sa\reg={\sp}\reg\otimes_k {\ss}\otimes_k kW$.

We denote by 
$$M\mapsto M\reg=\sa\reg\otimes_\sa M:\sa\mMOD\to \sa\reg\mMOD$$
the localisation functor. Note that
$\Res_{{\sp}\reg}M\reg={\sp}\reg\otimes_{\sp} M$. Note also that every
element of $M\reg$ can be written as 
$\alpha^r\otimes m$ for some $r\le 0$, $m\in M$, where
$\alpha=\prod_{H\in\CA}\alpha_H$. This makes the localisation functor
have specially good properties.

The restriction functor $\sa\reg\mMOD\to \sa\mMOD$ is a right adjoint to the
localisation functor. It is fully faithful. The adjunction morphism
coincides with the natural localisation
morphism $M\to M\reg$ of $\sa$-modules.
Its kernel is $M\tor$, the submodule of $M$ of
elements whose support is contained in $V-V\reg$.
Denote by $(\sa\mMOD)\tor$ the full subcategory of $\sa\mMOD$ of objects
$M$ such that $M\reg=0$. The following is clear.

\begin{lemma}
The localisation functor induces an equivalence
$$\sa\mMOD/(\sa\mMOD)\tor\iso \sa\reg\mMOD.$$
\end{lemma}

The category $\CO$ is a Serre subcategory of $\sa\mMOD$. 
Let $\CO\tor=\CO\cap (\sa\mMOD)\tor$. Then, the canonical functor
$\CO/\CO\tor\to \sa\mMOD/(\sa\mMOD)\tor$ is fully faithful.
Consequently,
the canonical functor
$$\CO/\CO\tor\to \sa\reg\mMOD$$
is fully faithful, with image a full abelian subcategory closed under
taking subobjects and quotients (but in general not closed under extensions).

\subsubsection{}
\label{diagduality}
When $k$ is a field,
we have a commutative diagram
$$\xymatrix{
D^b(\sa(\gamma)\mMod) \ar[rrrr]^{R\Hom_\sa(-,\sa)\ii}_\sim
\ar[d] &&&&
D^b(\sa(\gamma\ii)\mMod)^\opp \ar[d]\\
D^b(\sa(\gamma)\reg\mMod) \ar[rrrr]^{R\Hom_{\sa\reg}(-,\sa\reg)\ii}_\sim &&&&
D^b(\sa(\gamma\ii)\reg\mMod)^\opp
}$$
where the vertical arrows are given by localisation.

\subsubsection{}
\begin{lemma}
\label{dualitytorsion}
Assume $k$ is a field. Then, $(-)^\vee$ restricts to an equivalence
$\CO\tor(V,\gamma)\iso \CO\tor(V^*,\gamma\ii)^\opp$.
\end{lemma}

\begin{proof}
Let $M\in\CO$. We put a grading on $M$ (Proposition \ref{decomposition}).
Since $M$ is a finitely generated graded ${\sp}$-module (Lemma \ref{finiteness}),
the dimension of $\Ch(M)$, the characteristic variety of $M$, can
be obtained from the growth of the function $i\mapsto\dim M_i$. In particular,
$M\in\CO\tor$ if and only if
$\lim_{i\to\infty}\,\bigl(i^{1-\dim V}\cdot\dim M_i\bigr)=0$.
Such a property is preserved by $(-)^\vee$.
\end{proof}

Denote by $\gV: \CO\to\bar{\CO}=\CO/\CO\tor\,,\,
M\mapsto \Mbar,$ the quotient functor
(the notation $\gV$ has been used by Soergel \cite{So1}
for an analogous functor in the Lie algebra setup).
The  functor $\gV$ admits, 
by the standard `abstract nonsense', both a left adjoint and right
adjoint functors ${{}^\top{\gV}}, \gV^\top:\ \bar{\CO}
\to \CO$.

\begin{thm} 
\label{doublecentralizer}
Assume $k$ is a field, and  $Q$ is a projective in $\CO$.
Then, the canonical  adjunction morphism $a: Q\to \gV^\top(\Qbar)$ is
an isomorphism. In particular,  for any  object
$M$ in $\CO$, the  following canonical morphism  is an isomorphism
\begin{equation}\label{good}
\gV_*:\ \Hom_{\CO}(M,Q)\iso \Hom_{\bar{\CO}}(\Mbar,\Qbar).
\end{equation}
\end{thm}

\begin{proof} By  \S \ref{fullyfaithful}, for any  two objects
$M, Q,$ of $\CO$, we have a canonical isomorphism
$$\Hom_{\bar{\CO}}(\Mbar,\Qbar)\iso \Hom_{\sa\reg}(M\reg,Q\reg).$$

Assume $Q$ has a $\Delta$-filtration. Then it is free over ${\sp}$ and thus has no
non-zero submodule lying in $\CO\tor$, hence $\gV_*$ is injective.

Assume furthermore that $M$ has a $\nabla$-filtration. Then
$M^\vee$ has a $\Delta$-filtration
(Proposition \ref{usualduality}), hence has no non-zero submodule lying
in $\CO\tor$.
Since $(-)^\vee$ restricts to an equivalence
$\CO\tor(V,\gamma)\iso \CO\tor(V^*,\gamma\ii)^\opp$
(Lemma \ref{dualitytorsion}),
it follows that $M$ has no non-zero quotient lying in $\CO\tor$.
This shows that $\gV_*$ in \eqref{good} is an isomorphism.

From now on, we assume that  $Q$ is projective.
It follows that  $Q'=D(Q)$ is tilting
(Proposition \ref{Cher_dualityequivalences}), hence $\nabla$-filtered.

Now let $M$ be a  $\Delta$-filtered object. Then,
$M'=D(M)$ is $\Delta$-filtered. 
We apply the result  on $\gV_*$, that we have already proved, to
$\CO(\gamma\ii)$.
This yields,  by duality (cf \S \ref{diagduality}),
that \eqref{good} is an isomorphism,
 for any  $\Delta$-filtered object  $M$.

Since any projective is  $\Delta$-filtered,  for any
two projective objects $P,Q$ in $\CO$,
we have established the isomorphisms
\begin{equation}\label{good2}
\Hom_\CO(P, Q)\stackrel{\gV_*}\iso \Hom_{\bar{\CO}}(\Pbar,\Qbar)
\stackrel{\text{adjunction}}{=\!=\!=\!=\!=}
\Hom_\CO(P, \gV^\top(\Qbar)).
\end{equation}

The above isomorphisms imply, in particular, that, for any
indecomposable projective $P$,  we have
$\dim\Hom_\CO(P, Q)=\dim\Hom_\CO(P, \gV^\top(\Qbar))$.
It follows readily that
 the objects $Q$ and $\gV^\top(\Qbar)$
have the same composition factors with the same multiplicities.

We can finally prove that the canonical adjunction map
$a: Q\to \gV^\top(\Qbar)$ is an isomorphism.
By the previous paragraph, it suffices to show that $a$ is injective.
To this end, put $K:=\ker(a)$, and assume $K\neq 0$.
Let $L(E)$ be a simple submodule in $K$, and 
$P(E) \twoheadrightarrow
L(E),$ its projective cover.
By  construction,  the composite map
$g: P(E)\twoheadrightarrow
L(E)
\hookrightarrow K\hookrightarrow Q$ is nonzero.
This map $g \in \Hom_\CO(P(E) , Q)$ goes, under the
isomorphism between the left-hand and right-hand sides
of \eqref{good2}, to the map $a\ccirc g:\
P(E)\to K\hookrightarrow Q\stackrel{a}\to \gV^\top(\Qbar)$.
But the latter map is the zero map since $K=\ker(a)$,
which contradicts the fact that \eqref{good2} is an
isomorphism.
Thus, $\ker(a)=0$, and the Theorem is proved.
\end{proof}

\begin{rem}
In general, the assumption that $Q$ is projective cannot be replaced by
the weaker assumption that it is $\Delta$-filtered (already for
$W=\BZ/2\BZ$). Nevertheless, see 
Proposition \ref{opdam_add}.
\end{rem}

\begin{cor}\label{doublecentralizer_cor}
Let $X$ be a progenerator of $\CO$ and $\CE:=(\End_{\bar{\CO}} \Xbar)^\opp$.
Then there is
 an equivalence $\CO\iso~(\CE\mMod)^{\opp}.$
\end{cor}
\begin{proof} 
The preceding theorem implies that $(\End_{\CO} X)^\opp\iso\CE$ since
projective modules are $\Delta$-filtered. Hence we can use category 
equivalences \eqref{progen}.
\end{proof}

\pagebreak[3]
\pagebreak[3]\subsection{Dunkl operators}
\label{Dunkl}
\subsubsection{}
One has an $\sa$-action
 on the vector space ${\sp}$, hence an $\ss$-action, arising via
the identification $\sp=\Delta(k)$.
One finds, in particular, that
the action of $\xi\in V$ on ${\sp}$ is given by the
Dunkl operator
$$T_\xi=\partial_\xi+\sum_{H\in\CA}\frac{\langle\xi,
\alpha_H\rangle}{\alpha_H}{\cdot} a_H \;\in\;\CD(V\reg)\rtimes W,$$
where $\CD(V\reg)$ stands for the algebra of regular differential
operators on $V\reg$, acted upon by $W$ in a natural way,
and
 $a_H\in kW$ is viewed as an element of $\CD(V\reg)\rtimes W$.
It follows that $T_\xi(\sp)\subset \sp$ (as part of $\sa$-action on
$\sp=\Delta(k)$); furthermore,
this $\sa$-action on $\sp$ is known (Cherednik,
\cite[Proposition 4.5]{EtGi})
to be faithful:

\begin{thm}
\label{diffop}
The $\sa$-representation $\Delta(k)$ is faithful. Thus,
the natural action of ${\sp}W$ on ${\sp}$ extends to
an injective algebra morphism
$i:\sa\hookrightarrow k\otimes_\BC \CD(V\reg)\rtimes W$ which
maps $\xi\in V$ to $T_\xi$.

The map $i$  induces an algebra  isomorphism
$\sa\reg\iso k\otimes_\BC \CD(V\reg)\rtimes W$.
\end{thm}

\subsubsection{}

We consider $M=\Ind_{{\ss}W}^\sa X={\sp}\otimes X$,
where $X$ is locally nilpotent and finitely generated as an
${\ss}$-module, free over $k$.
The action of $\xi\in V$ on $p\otimes v$, $p\in {\sp}$ and $v\in X$ is given by
$$\xi(p\otimes v)=p\otimes\xi v+\partial_\xi (p)\otimes v+
\sum_H \sum_{0\le i,j\le e_H-1} e_H(k_{H,i+j}-k_{H,j})
\frac{\alpha_H(\xi)}{\alpha_H} \varepsilon_{H,i}(p)\otimes \varepsilon_{H,j}(v).$$
Using Dunkl operators, \ie, via the isomorphism of Theorem \ref{diffop},
we have a structure of $W$-equivariant $(k\otimes_\BC\CD(V\reg))$-module
on $M\reg$.
The corresponding connection is given by
$$\partial_\xi=\xi-\sum_H \frac{\alpha_H(\xi)}{\alpha_H}\cdot
\bigl(\sum_{0\le i,j\le
e_H-1} e_H k_{H,i+j}\eps_{H,i}\otimes \eps_{H,j}\bigr).$$
Hence
$$\partial_\xi(p\otimes v)=
p\otimes \xi v+\partial_\xi(p)\otimes v-\sum_H \sum_i
e_H k_{H,i}\-\frac{\alpha_H(\xi)}{\alpha_H}p\otimes \eps_{H,i}(v).$$

For the rest of \S \ref{sectionHecke}, we assume that $k=\BC$.

The following result is well-known to experts, but we could not find
an appropriate reference in the literature.
\begin{prop}
The above formula for $\partial_\xi$ defines a $W$-equivariant integrable algebraic
connection on $M$ with regular singularities.
\end{prop}
\begin{proof}
All the claims follow  from the construction, with the
exception of the assertion that the singularities of the
connection are regular. The connection has visibly only simple
poles at the reflection hyperplanes, hence it suffices to prove
the regularity at infinity with respect to some (hence any, see
[De]) compactification of $V$.

Consider the  $W$-equivariant compactification
$Y=\mathbb{P}(\mathbf{C}+V)$ of $V$, and extend $M$ to the free
$\mathcal{O}_Y$-module $M_Y:=\mathcal{O}_Y\otimes X$.
Using a filtration of $X$ we can reduce to the
case where $X=E$ is simple.

A straightforward computation shows
that with respect to the extension $M_Y$ of $M$ and with respect
to any standard coordinate patch on $Y$, the poles at infinity are
also simple in this case.
\end{proof}

\subsubsection{}
\label{subsubmonodromy}
We define a morphism of abelian groups $r:K_0(\CO)\to\BZ$
by $r([\Delta(E)])=\dim E,$ for $E\in\Irr(\BC W)$.

\begin{lemma}
\label{rankvectorbundle}
Let $M\in\CO$. Then, $M\reg$ is a vector bundle of rank $r([M])$
on $V\reg$.
\end{lemma}

\begin{proof}
Since $M\reg$ is a
finitely generated $\BC[V\reg]$-module with a connection, it is a 
vector bundle.

Now, taking the rank of that vector bundle induces a morphism
$K_0(\CO)\to\BZ$, which takes the correct value on $\Delta(E)$.
\end{proof}

\subsubsection{}
\begin{prop}\label{opdam_add} 
Assume  $k_{H,i}-k_{H,j}+\frac{i-j}{e_H}\not\in\BZ$,
for all $H\in\CA$ and all $0\le i\not= j\le e_{H}-1$.
Let $N$ be a $\Delta$-filtered object in $\CO$.
Then, for any $M\in\CO$, we have
$\Hom_\CO(M,N)\iso \Hom_{\bar{\CO}}(\bar{M},\bar{N})$.
\end{prop}

\begin{proof} Assume first that $M$ is also 
 a $\Delta$-filtered object.
Then, we can write $M=\operatorname{Ind}_{{\ss}W}^{\sa} X$ and
$N=\operatorname{Ind}_{{\ss}W}^{\sa} Y$ with $X,Y$ finite dimensional
${\ss}W$-modules, nilpotent over ${\ss}$. The space $\Hom_{\sa}(M,N)$ is the
intersection of $\Hom_{\sp}(M,N)={\sp}\otimes\Hom_k(X,Y)$ with
$\Hom_{\sa\reg}(M\reg,N\reg)$. As in the proof of Theorem
\ref{doublecentralizer}, we have to show that any element $\Psi$
of $\Hom_{{\sp}\reg}(M\reg,N\reg)$ that commutes with the action
of $\sa\reg$ extends to a $\sp$-morphism $M\to N$. Observe that
$\Psi$ is nothing but a flat, $W$-invariant section of the
connection on $\Hom_{\sa\reg}(M\reg,N\reg)$.

The residue of this connection on a hyperplane $H\in\CA$ is
constant, and has eigenvalue $e_H(k_{H,i}-k_{H,j})$ on
$\Hom_k(X_i,Y_j)$, where $X_i$ is the summand of $X$ of $W_H$-type
$\operatorname{det}_{|W_H}^{-i}$ (and likewise for $Y_j$).

Locally near a generic point $p$ of $H$ we expand $\Psi=
\sum_{l\geq l_0}\alpha_H^l\Psi_l$ with $\Psi_l$ holomorphic on $H$
near $p$, of $W_H$-type $\operatorname{det}_{|W_H}^{l}$, and with
$\Psi_{l_0}$ not identically zero on $H$. From the lowest order
term of the equation $\partial_{v_H}(\Psi)=0$ we see that there
exist $i,j$ such that $i-j =
l_0\operatorname{mod}(e_H\BZ)$, and such that
$l_0+e_H(k_{H,i}-k_{H,j})=0$. Thus $l_0=0$ and $\Psi$ is regular
on $H$.
This completes the proof of the Proposition in the special
case where both $M$ and $N$ are $\Delta$-filtered.

The general case follows from the special case above by repeating the
part of the argument from the
 proof of Theorem \ref{doublecentralizer}, starting with formula
\eqref{good2}.
\end{proof}
\begin{rem}
The condition of the Proposition is equivalent to the semi-simplicity of
the Hecke algebra $\CH(W_H)$ of $W_H$. One could conjecture that this
assumption can be replaced by the assumption that 
$\kz(N)$ is a projective  $\CH(W_H)$-module (this would still not cover completely
Theorem \ref{doublecentralizer}).
\end{rem}
\begin{rem}
If $e_H=2$ for all $H$, then the condition of the Proposition reads:
$k_H\not\in \frac{1}{2}+\BZ$.
\end{rem}

\subsubsection{}
\label{definitionHecke}
Let $\BC[(\Bk_{H,i})_{1\le i\le e_H-1}]$ be the polynomial ring in the
indeterminates $\Bk_{H,i}$ with $\Bk_{w(H),i}=\Bk_{H,i}$ for $w\in W$.
We have a canonical morphism of $\BC$-algebras
$\BC[(\Bk_{H,i})]\to \BC, \ \Bk_{H,i}\mapsto k_{H,i}$.
Let $\Gm$ be the kernel of that morphism and $R$ the completion of
$\BC[(\Bk_{H,i})]$ at the maximal ideal $\Gm$.

Fix $x_0\in V\reg$, and let $B_W=\pi_1(V\reg/W,x_0)$ be the Artin braid
group associated to $W$.

Let $\CH_R=\CH_R(W,V,\gamma)$ be the Hecke algebra of $W$ over $R$,
that is the quotient of
$R[B_W]$ by the relations
$$(T-1)\prod_{j=1}^{e_H-1} (T-\det(s)^{-j}\cdot e^{2i\pi \Bk_{H,j}})=0$$
for $H\in\CA$, $s\in W$ the reflection around $H$ with non-trivial eigenvalue
$e^{2i\pi/e_H}$ and $T$ an $s$-generator of the monodromy around $H$,
cf \cite[\S 4.C]{BrMaRou}. Note that the parameters differ from
\cite{BrMaRou} because
we will be using the horizontal sections functor instead of the solution
functor.

We put $\CH_K=\CH_R\otimes_R K$, where $K$ is the field of
fractions of $R$ and $\CH=\CH_R\otimes_R (R/\Gm)$.

\begin{rem}\label{W}
It is known that $\CH_R$ is free of rank $|W|$ over $R$ for
all $W$ that do not have an irreducible component of type
$G_{17\ldots 19}$, $G_{24\ldots 27}$,
$G_{29}$, $G_{31\ldots 34}$ in Shephard-Todd notation
(in these cases, the statement is conjectural) \cite{Mu}.
\end{rem}

\subsection{The Knizhnik-Zamolodchikov functor.}
Let $M$ be a $(\BC[V\reg]\rtimes W)$-module, free of finite rank over
${\sp}\reg=\BC[V\reg]$.
Let $\nabla:M\to M\otimes_{\BC} R$ be an $R$-linear integrable connection.
Then, the horizontal sections of $\nabla$ define,
via the monodromy representation, an
$RB_W$-module $L$, free over $R$.

Let $\nabla_0:M\to M$ be the special fiber of $\nabla$. Then, the
horizontal sections of $\nabla_0$ is the $\BC B_W$-module
$L\otimes_R (R/\Gm)$.

Let $\nabla_K:M\to K\otimes_\BC M$ be the generic fiber of $\nabla$.
Then, the horizontal sections of $\nabla_K$ is the $KB_W$-module
$L\otimes_R K$.

\smallskip
Taking horizontal sections defines an exact functor from the category 
of $W$-equivariant vector bundles on $R\otimes_\BC V\reg$ with an integrable
connection to the category of $RB_W$-modules that are
free of finite rank over $R$.

Since the connection on $\Delta(R\otimes_\BC E)\reg$ has regular
singularities it follows that the connection on $M\reg$
has regular singularities for any $M\in\CO_R^\Delta$.

Composing with the localisation functor, we obtain an exact
functor ${\kz}_R$ from  $\CO_R^\Delta$ to
the category of $RB_W$-modules that are
free of finite rank over $R$.

Similarly, we obtain exact functors
${\kz}:\CO\to \BC B_W\mMod$ and 
${\kz}_K:\CO_K\to KB_W\mMod$.

It is well-known (cf. e.g. \cite[Theorem 4.12]{BrMaRou}) that
the representation of $KB_W$ on ${\kz}_K(\Delta(K\otimes_\BC E))$ 
factors through $\CH_K$ to give a representation corresponding (via
Tits' deformation Theorem) to the representation $E$ of $\BC W$.
Recall that $\CH=\CH_R\otimes_R (R/\Gm)$.

\begin{thm}[{\bf Hecke algebra action}]
\label{factorHecke}
The functor ${\kz}:\CO\to \BC B_W\mMod$ factors through a functor
${\kz}:\CO/\CO\tor\to \CH\mMod$.
Similarly, the functor ${\kz}_K:\CO_K\to KB_W\mMod$ factors through a functor
${\kz}_K:\CO_K/(\CO_K)\tor\to \CH_K\mMod$.

For $M\in\CO_R^\Delta$, the action of
$RB_W$ on ${\kz}_R(M)$ factors through $\CH_R$.

We have a commutative diagram
$$\xymatrix{
\CO_K \ar[rr]^{{\kz}_K} && \CH_K\mMod  \\
\CO_R^\Delta \ar[rr]^{{\kz}_R}\ar[d]_{\BC\otimes_R -}\ar[u]^{K\otimes_R -} &&
 \CH_R\mMod \ar[d]^{\BC\otimes_R -}\ar[u]_{K\otimes_R -} \\
\CO \ar[rr]^{{\kz}} && \CH\mMod 
}$$
\end{thm}

\begin{proof}
First, $\CO\tor$ (and $(\CO_K)\tor$) are the kernels of localisation.

When $M=\Delta(K\otimes_\BC E)$, then, we have the Knizhnik-Zamolodchikov
connection and the representation ${\kz}_K(M)$ factors through $\CH_K$.
Since $\CO_K$ is semi-simple (Corollary \ref{semisimple}), it
follows that the action on
${\kz}_K(M)$ factors through $\CH_K$ for any $M$ in $\CO_K$.

We now consider the case of a $\Delta$-filtered module $M$ of $\CO_R$.
We know that the action of $KB_W$ on
$K\otimes_R {\kz}_R(M)\simeq {\kz}_K(K\otimes_R M)$
factors through $\CH_K$.
Since ${\kz}_R(M)$ is free over $R$, it follows that the action
of $RB_W$ on ${\kz}_R(M)$ factors through $\CH_R$.

   From this result, we deduce that the action of $\BC B_W$ on
${\kz}(\Delta^r(\BC W))\iso \BC\otimes_R {\kz}_R(\Delta^r(RW))$ factors through
$\CH$. Since every indecomposable projective object of $\CO$ is a
direct summand of $\Delta^r(\BC W)$ for appropriate $r$
(Corollary \ref{projective}), it follows
that the action of $\BC B_W$ on ${\kz}(M)$ factors through $\CH$ for
every projective $M$, hence for every $M$ in $\CO$.
\end{proof}

\pagebreak[3]\subsection{Main results}\label{main}
In this subsection we assume  that  $\dim \CH=|W|$, cf. Remark \ref{W}.

The functor  ${\kz}:\CO\to \CH\mMod$ is exact.
Hence, it is represented by a projective $P_\kz\in \CO$.
In other words, there exists an algebra 
 morphism $\phi:\CH\to(\End_\CO\,P_\kz)^\opp$ such that
the functor ${\kz}$ is isomorphic to $\Hom_\CO(P_\kz,-)$.

We know also, see \S \ref{fullyfaithful}, that the functor
${\kz}$ factors through $\CO/\CO\tor\to \CH\mMod$.

\begin{thm}
\label{OfromH}
The functor ${\kz}$ induces an  equivalence: $\CO/\CO\tor\iso \CH\mMod$.
\end{thm} 

This Theorem  is equivalent
to
\begin{thm}
\label{OfromH2} The morphism  $\phi:\CH\to(\End_\CO\,P_\kz)^\opp$
is an algebra isomorphism.
\end{thm}
\begin{proof}[Proof of Theorems \ref{OfromH}-\ref{OfromH2}.]
Recall that the horizontal sections functor gives an equivalence
from the category of vector bundles over $V\reg/W$ with
a regular integrable connection to the category of finite-dimensional
$\BC B_W$-modules
(Riemann-Hilbert correspondence, \cite[Theorems I.2.17 and II.5.9]{De}).

We deduce from \S \ref{fullyfaithful} that 
${\kz}:\CO/\CO\tor\to \CH$ is a fully faithful exact functor with image
a full subcategory closed under taking subobjects and quotients.
Furthermore, ${\overline{P_\kz}}$, the image of
$P_\kz$ in $\CO/\CO\tor$,  is a progenerator  of $\CO/\CO\tor$.
Thus, Theorem \ref{OfromH} follows from Theorem \ref{OfromH2}.

To prove Theorem \ref{OfromH2},
observe  that the morphism $\phi$ is surjective. Indeed,
let $\CC'$ be a full subcategory of an abelian category $\CC$,
closed under taking quotients, and ${{}^{\top\!}F}$ a left adjoint to the inclusion
$F:
\CC'\hookrightarrow \CC$.
Then the adjunction morphism $\eta:\Id_\CC\to F\ccirc ({{}^{\top\!}F})$
 is surjective, since $\eta(X)$ is the
canonical map from $X$ to its largest quotient in $\CC'$.
This proves surjectivity of the morphism $\phi$ above.

Further, we have
$$P_\kz=\bigoplus_{E\in\Irr(\BC W)} (\dim {\kz} (L(E))) P(E).$$
Hence, we compute
$$\dim(\End_\CO\,P_\kz)=\bigoplus_{E,F} \dim {\kz}(L(E)) \dim {\kz}(L(F))\dim \Hom(P(E),P(F))
$$
$$=\bigoplus_{E,F,G}\dim {\kz}(L(E)) \dim {\kz}(L(F))[P(E):\Delta(G)][\Delta(G):L(F)]$$
$$=\bigoplus_{E,F,G}\dim {\kz}(L(E)) \dim {\kz}(L(F))[\nabla(G):L(E)][\Delta(G):L(F)]$$
$$=\bigoplus_G \dim {\kz}(\nabla(G)) \dim {\kz}(\Delta(G))$$
Now, the restrictions of $\nabla(G)$ and $\Delta(G)$
to $V\reg$ are vector bundles of rank $\dim G$ (Proposition
\ref{classdeltanabla} and Lemma \ref{rankvectorbundle}), hence
$\dim(\End_\CO\,P_\kz)=|W|=\dim \CH$.
This shows that $\phi$ is an isomorphism.
Note that this rank computation can also be achieved by deformation to
$R$.
\end{proof}
The following result shows that the category
$\CO$ can be completely recovered from $\CH$
and a certain $\CH$-module~:

\begin{thm}[{\bf Double-centralizer property}]
\label{OfromH3} Let $Q$ be a projective in $\CO$. Then,  the  canonical map
$\Hom_\CO(M,Q)\to \Hom_{\CH}\bigl(\kz(M)\,,\,\kz(Q)\bigr)$ is an
isomorphism, for any $M\in \CO$.

 Furthermore, if $X$ is a progenerator for $\CO$, then,
we have an equivalence
$$\bigl(\End_{\CH}\,{\kz}(X)\bigr)^\opp \mMod\iso \CO.$$
\end{thm}
\begin{proof} The first part follows from Theorems \ref{doublecentralizer} and \ref{OfromH2}
and the second from Corollary  \ref{doublecentralizer_cor}.
\end{proof}

\begin{rem}\label{qSchur} We conjecture that, if $W=S_n$, then
$\CO$ is equivalent to the category of finitely-generated modules over the
associated $q$-Schur algebra. That would imply, in particular,
that if $\,k_{H,1}=k_1<0\,$ is a negative real constant, then
the  Cherednik algebras $\sa(S_n)$ with parameters  $k_1$ and
$k_1-1,$ respectively, are Morita equivalent.
\end{rem}

Let $Z(\CH)$ denote the center of the algebra $\CH$ and
$Z(\CO)$  the {\it center} of category $\CO$ (i.e. the
algebra of endomorphisms of the identity functor $\Id_\CO$).
\begin{cor}
The canonical morphism $Z(\CO)\to\End_{\CO} P_{\kz}$ induces an
isomorphism $Z(\CO)\iso Z(\CH)$.
In particular, the functor $\kz$ induces a bijection between blocks of
$\CO$ and blocks of $\CH$.
\end{cor}

\begin{proof}
This follows immediately from Theorem \ref{OfromH3}~: given two
rings $B$ and $C$ and a $(B,C)$-bimodule $M$ such that the canonical
morphisms $B\iso \End_{C^\opp}(M)$ and $C\iso (\End_B M)^\opp$ are
isomorphisms, then we have a canonical isomorphism
$Z(B)\iso Z(C)$.
\end{proof}

The decomposition matrix $K_0(\CO_K)\to K_0(\CO)$
is triangular. We deduce the triangularity of decomposition matrices
of Hecke algebras, in characteristic $0$~:

\begin{cor}
\label{decHecke}
The decomposition matrix $K_0(\CH_K)\to K_0(\CH)$ is triangular.
\end{cor}

\subsubsection{$\kz$-functor and Twist}
\label{twist}
Let $\zeta$ be a one-dimensional character of $W$ and 
$\tau_{_\zeta}:\BC W\iso \BC W$ the automorphism given by
$w\mapsto \zeta(w)\cdot w$ for $w\in W$.
This extends to an isomorphism
$$\tau_{_\zeta}:A(\gamma)\iso A(\tau_{_\zeta}(\gamma)),\ \
V\ni\xi\mapsto \xi,\ \ V^*\ni x\mapsto x,\ \  W\ni w\mapsto \zeta(w)\cdot w.$$
We obtain an equivalence
$\CO(\gamma)\iso \CO(\tau_{_\zeta}(\gamma))$, sending 
$V(E)$ to $V(E\otimes \zeta^{-1})$, where $V$ stands for any of the symbols:
$L,\Delta,\nabla,P,I,T$.

For $H\in\CA$, let $d_H\in\{1,\ldots,e_H\}$ such that
$\zeta_{|W_H}=\det_{|W_H}^{d_H}$.
Define an automorphism $\eta_{_\zeta}$ of $\CD(V_{\reg})\rtimes W$
by
$$P\ni f\mapsto f, \ \ W\ni w\mapsto \zeta(w)\cdot w$$
$$\textrm{and }\partial_\xi\mapsto \partial_\xi-
\sum_H \frac{\alpha_H(\xi)}{\alpha_H}\eps_{H,e_H}\cdot e_H\cdot  k_{H,e_H-d_H}
\textrm{ for }\xi\in V.$$
(for notation, see Remark \ref{eps}).
We have a commutative diagram
$$\xymatrix{
A(\gamma)\ar[rr]^i \ar[d]_{\tau_{_\zeta}}^\sim &&
 \CD(V_{\reg})\rtimes W \ar[d]^{\eta_{_\zeta}}_\sim \\
A(\tau_{_\zeta}(\gamma))\ar[rr]_i && \CD(V_{\reg})\rtimes W
}$$

Given $M$ a $(\CD(V_{\reg})\rtimes W)$-module, then
$(\eta_{_\zeta})_* M\iso M\otimes_{\CO_{V_{\reg}}} \Delta(\zeta^{-1})_{\reg}$.

This self-equivalence of the category of $W$-equivariant bundles with
a regular singular connection on $V_{reg}$ corresponds, via the
horizontal sections functor, to the automorphism of
$\BC B_W$ given by 
$$T\mapsto e^{-2i\pi k_{H,e_H-d_H}}\zeta(s)^{-1} T$$
for $H\in\CA$, $s\in W$ the reflection around $H$ with non-trivial eigenvalue
$e^{2i\pi/e_H}$ and $T$ an $s$-generator of the monodromy around $H$.
This induces an isomorphism
$\CH(\zeta):\CH(W,\gamma)\iso \CH(W,\tau_{_\zeta}(\gamma))$ and
the following diagram is commutative~:
$$\xymatrix{
\CO(\gamma)\ar[rr]^{(\tau_{{_\zeta}})_*}_\sim \ar[d]_{\kz} &&
 \CO(\tau_{_\zeta}(\gamma))\ar[d]^{\kz} \\
\CH(W,\gamma)\mMod\ar[rr]_{\CH(\zeta)}^\sim && \CH(W,\tau_{_\zeta}(\gamma))\mMod
}$$

\subsubsection{$\kz$-functor and Duality}
\label{KZduality}
We have a commutative diagram
$$\xymatrix{
D^b(\sa(\gamma)\mMod) \ar[rrrr]^{(\tau_{\det})_*\circ R\Hom_{\sa(\gamma)}(-,\sa(\gamma))^\dag}_\sim
\ar[d] &&&&
D^b(\sa(\tau_{\det}(\gamma\ii))\mMod)^\opp \ar[d]\\
D^b((\CD(V\reg)\rtimes W)\mcoh) \ar[rrrr]^\sim &&&&
D^b((\CD(V\reg)\rtimes W)\mcoh)^\opp
}$$
where the vertical arrows are given by localisation
followed by the Dunkl operator isomorphism $i$ of Theorem \ref{diffop} and
the bottom horizontal arrow is the classical $\CD$-module duality.

Consider the isomorphism 
$\BC B_W\iso (\BC B_W)^\opp$ given by
$T\mapsto \det(s)^{-1}e^{2i\pi k_{H,1}}T^{-1}$

for $H\in\CA$, $s\in W$ the reflection around $H$ with non-trivial eigenvalue
$e^{2i\pi/e_H}$ and $T$ an $s$-generator of the monodromy around $H$.
It induces an isomorphism
$$\CH({}\ii):\;\CH(W,\gamma)\iso \CH(W,\gamma\ii)^\opp.$$

We conclude that we have a commutative diagram
$$\xymatrix{
D^b(\CO(\gamma))\ar[rrrr]_\sim^D \ar[d]_{{\kz}} &&&&
D^b(\CO(\gamma\ii))^\opp \ar[d]^{{\kz}} \\
D^b(\CH(W,\gamma))\ar[rrrr]^\sim_{\CH(\ii)} &&&& D^b(\CH(W,\gamma\ii))^\opp
}$$

\smallskip
On the other hand, by Lemma \ref{dualitytorsion}, we know that
$(-)^\vee$ preserves
$\CO\tor$, hence descends to the quotient category $\CO/\CO\tor$, i.e.,
 there is an equivalence $\Phi$ making the
following diagram commute:
$$\xymatrix{
\CO(V,\gamma)\ar[rrr]_\sim^{-^\vee} \ar[d]_{{\kz}} &&&
\CO(V^*,\gamma\ii)^\opp \ar[d]^{{\kz}} \\
\CH(W,V,\gamma)\mMod\ar[rrr]^\sim_{\Phi} &&& \CH(W,V^*,\gamma\ii)\mMod^\opp
}$$

Further, choose a $W$-invariant hermitian form on $V$, i.e.,
a semi-linear $W$-equivariant isomorphism $\kappa:V\iso V^*$.
Then, we get an isomorphism
$\pi_1(V\reg/W,x_0)\iso \pi_1(V^*\reg/W,\kappa(x_0))$.
It induces an isomorphism
$\CH(\kappa):\CH(W,V,\gamma)\iso \CH(W,V^*,\gamma)$. Composing
with $\CH({}\ii)$, we obtain an isomorphism
$\CH(\kappa\ccirc(-)\ii):\CH(W,V,\gamma)\iso \CH(W,V^*,\gamma\ii)^\opp$,
which we denote below by $\psi$.

\begin{rem} One could conjecture that the two functors
$\Phi$ and $\psi_*$ are isomorphic
(they induce the same maps at the level of Grothendieck groups).
\end{rem}
\subsubsection{The $\sa$-module $P_{\kz}$ and Duality}
\label{P_KZ}
Let $\Irr(W,V,\gamma)\subset \Irr(W)$ denote  the subset
formed by all $E\in \Irr(W)$ such that $L(E)\reg\not=0$.
We have a bijection $\Irr(W,V,\gamma)\iso 
\Irr(W,V^*,\gamma\ii)\,,\, E\mapsto E^\vee$
(Proposition \ref{usualduality} and Lemma \ref{dualitytorsion}).
Thus, $P_{\kz}=\bigoplus_{E\in \Irr(W,V,\gamma)} (\dim \kz(L(E)))\cdot P(E)$.

To make the dependence on $V$ and $\gamma$ explicit,
we will write $P_{\kz}=P_{\kz}(V,\gamma)$.

\begin{prop} {\sf{(i)}}\; We have $D(P_{\kz}(V,\gamma))\simeq
P_{\kz}(V,\gamma\ii)$
and  $P_{\kz}(V,\gamma)^\vee \simeq P_{\kz}(V^*,\gamma)$. In particular,
 $P_{\kz}$ is both projective and injective.

{\sf{(ii)}}\; For $E\in\Irr(W)$, the following  are equivalent
\begin{itemize}
\item $E\in \Irr(W,V,\gamma)$
\item $L(E)$ is a submodule of a standard module
\item $P(E)$ is a submodule of $P_{\kz}$
\item $P(E)$ is injective
\item $P(E)$ is tilting
\item $I(E)$ is projective
\item $I(E)$ is tilting 
\end{itemize}
\end{prop}

\begin{proof} The first claim
follows from \S \ref{KZduality}.
 Proposition \ref{usualduality} then implies that $P_{\kz}$ is injective.

The considerations above  imply that if $E\in \Irr(W,V,\gamma)$, then
$P(E)$ is injective and tilting. The assertions about
$I(E)$ follow by applying $(-)^\vee$.

We know that if $L(E)$ is a submodule of a $\Delta$-filtered module
or a quotient of a $\nabla$-filtered module, then $E\in \Irr(W,V,\gamma).$

This shows that any of the assertions about $P(E)$ or $I(E)$ implies
that $E\in \Irr(W,V,\gamma).$
\end{proof}

\section{Relation to Kazhdan-Lusztig theory of cells}
\label{subKL}
We review some  parts of Kazhdan-Lusztig and Lusztig's theory of
Weyl group representations.

\pagebreak[3]\subsection{Lusztig's algebra $\CJ$}

\subsubsection{}
Let $(W,S)$ be a finite Weyl group, $\CH$ be its Hecke algebra, a
 $\BZ[v,v^{-1}]$-algebra with basis $\{T_w\}_{w\in W}$
and relations
$$T_wT_{w'}=T_{ww'} \text{ if }l(ww')=l(w)+l(w')\quad\text{and}\quad
(T_s+1)(T_s-v^2)=0 \text{ for }s\in S.$$

Lusztig associated to $W$ a  $\BZ$-ring $\CJ$,
usually referred to as {\it asymptotic Hecke algebra},
 \cite[\S 2.3]{Lu3}.
Let
$\varpi:\CH\to \BZ[v,v^{-1}]\otimes_\BZ \CJ$ be Lusztig's morphism of
$\BZ[v,v^{-1}]$-algebras
\cite[\S 2.4]{Lu3}.

The ring $\BQ\otimes_\BZ \CJ$ is semi-simple and
the morphism $\Id_{\BQ(v)}\otimes \varpi$ is an isomorphism. 

For any commutative $\BQ[v,v^{-1}]$-algebra $R$ we put $\CH_R:= R\otimes_{\BZ[v,v^{-1}]}\CH$.

\begin{defi}
The $\CH_R$-modules $S(M)=\varpi^*(R\otimes_\BQ M)$,
for $M$ a simple $\BQ\otimes_\BZ \CJ$-module,
will be referred to as {\sl standard $\CH_R$-modules}.\q{There seems to
be no name for such modules in the literature.}
\end{defi}

When $R=\BQ(v)$, then the standard $\CH_R$-modules are simple and
this gives a bijection from the set of simple $(\BQ\otimes_\BZ \CJ)$-modules
to the set of simple $(\BQ(v)\otimes_{\BZ[v,v^{-1}]}\CH)$-modules.

Similarly, taking $K=\BQ[v,v^{-1}]/(v-1)$, we obtain a bijection
from the set of simple $(\BQ\otimes_\BZ \CJ)$-modules
to the set of simple $\BQ W$-modules.

We will identify these sets of simple modules via these bijections.

We have an order $\le_{LR}$ on $W$ constructed in \cite[p.167]{KaLu}.
We denote by $\CC$ the set of two-sided cells of $W$ and by 
$\le$ the order on $\CC$ coming from $\le_{LR}$.

Let $\{C_w\}_{w\in W}$ be the Kazhdan-Lusztig basis for $\CH$.
Let $I$ be an ideal of $\CC$, \ie, a subset such that given $\uc\le\uc'$,
then $\uc'\in I\Rightarrow \uc\in I$.
We put
$\CH^I=\bigoplus_{\uc\in I,w\in\uc} \BZ[v,v^{-1}]C_w$.
This is a two-sided ideal of $\CH$ \cite[p.137]{Lu1}.

The ring $\CJ$ comes with a $\BZ$-basis $\{t_w\}_{w\in W}$ and we put
$\CJ_{\uc}=\bigoplus_{w\in\uc}\BZ t_w$. This is a block of $\CJ$ and
$\CJ=\bigoplus_{\uc\in\CC} \CJ_{\uc}$. The corresponding partition
of the set of simple $(\BQ\otimes_\BZ \CJ)$-modules is called the
partition into families.

Given $I$ an ideal of $\CC$, we denote by
$I^\circ$ the set of $\uc\in I$ such that there is $\uc'\in I$ with
$\uc<\uc'$.

The following is a slight reformulation of \cite[\S 1.4]{Lu3}~:

\begin{prop}
\label{filtrationJ}
Let $I$ be an ideal of $\CC$.
Then, the assignment $t_w\mapsto C_w$ induces an isomorphism of
$\CH$-modules
$$\bigoplus_{\uc\in
I-I^\circ}\varpi^*\bigl(\BZ[v,v^{-1}]\otimes_{\BZ}\CJ_{\uc}\bigr)
\iso
\CH^I/\CH^{I^\circ}.$$

In particular, the $(\BQ[v,v^{-1}]\otimes_{\BZ[v,v^{-1}]} \CH)$-module 
$\BQ[v,v^{-1}]\otimes_{\BZ[v,v^{-1}]} (\CH^I/\CH^{I^\circ})$
is a direct sum of standard $\CH_{\BQ[v,v^{-1}]}$-modules.
\end{prop}

This proposition gives a characterization of standard $\CH_{\BQ[v,v^{-1}]}$-modules
via the Hecke algebra filtration coming from two-sided cells.

\subsubsection{}
Next, we consider  filtrations coming from certain functions on the set
of two-sided cells.

\begin{defi}
A {\em sorting} function 
$f:W\to\BZ$ is a function constant on two-sided cells and such that 
$\uc'<\uc\Rightarrow f(\uc')>f(\uc)$.
\end{defi}

Given a sorting function $f$, we
put  $\CH^{\ge i}_R:=\bigoplus_{w\in W, f(w)\ge i} R\cdot C_w$ and
 $\CH^{>i}_R:=\bigoplus_{w\in W, f(w)>i} {R\cdot C_w.}$
Then, $\CH^{\ge i}_R$ is a two-sided ideal of $\CH_R$, since
$I=\{\uc\in\CC \;|\; f(\uc)\ge i\}$ is an ideal. Similarly,
$\CH^{>i}_R$ is a two-sided ideal of $\CH_R$.
Furthermore,
$I^\circ\subseteq \{\uc\in\CC \;|\; f(\uc)>i\}$.
Consequently, we deduce from Proposition \ref{filtrationJ}~:

\begin{cor}
We have an isomorphism of
$\CH$-modules
$$\bigoplus_{\uc\in\CC,f(\uc)=i}\varpi^*\bigl(\BZ[v,v^{-1}]\otimes_{\BZ}\CJ_{\uc}\bigr)\iso
\CH^{\ge i}/\CH^{>i}.$$

In particular, the $(\BQ[v,v^{-1}]\otimes_{\BZ[v,v^{-1}]} \CH)$-module 
$\BQ[v,v^{-1}]\otimes_{\BZ[v,v^{-1}]} (\CH^{\ge i}/\CH^{>i})$
is a direct sum of standard $\CH_{\BQ[v,v^{-1}]}$-modules.
\end{cor}

Thus, we have another characterization of standard $\CH_{\BQ[v,v^{-1}]}$-modules via
the Hecke algebra filtration coming from $f$.

\smallskip
Let $\CF$ be the set of families of irreducible characters of $W$.
We transfer the concepts associated with $\CC$ to $\CF$ via the canonical
bijection between $\CC$ and $\CF$.

In particular, we have a function $f:\Irr(W)\to\BZ$ constant on
families.

We have $\CH^{\ge i}=\CH\cap (\bigoplus_{f(E)\ge i} 
e_E \BQ(v)\otimes_{\BZ[v,v^{-1}]} \CH)$,
where $e_E$ is the primitive central idempotent of
$\BQ(v)\otimes_{\BZ[v,v^{-1}]} \CH$
that acts as $1$ on the simple $(\BQ(v)\otimes_{\BZ[v,v^{-1}]} \CH)$-module
corresponding to $E$.

This shows that, if $R$ is a localisation of $\BQ[v,v^{-1}]$, then
the filtration on $\CH_R=R\otimes_{\BZ[v,v^{-1}]} \CH$ given by $f$ can be recovered
without using the Kazhdan-Lusztig basis.
We obtain

\begin{prop}
\label{filtrationprojHecke}
Let $R$ be  a localisation of $\BQ[v,v^{-1}]$ and $P$ be a projective
$\CH_R$-module.
Let $Q^{\ge i}$ (resp. $Q^{>i}$) be the sum
of the simple submodules $E$ of $\BQ(v)\otimes_R P$ such that
$f(E)\ge i$ (resp. $f(E)>i$).

Then, $(P\cap Q^{\ge i})/(P\cap Q^{>i})$ is a direct sum of
standard $\CH_R$-modules.\hfill$\Box$
\end{prop}

Thus, any sorting function yields a characterization of the
standard $\CH_R$-modules without using the Kazhdan-Lusztig basis.

\subsubsection{}
Given $E\in\Irr(W)$, we denote by $a_E$ (resp. $A_E$) the lowest
(resp. highest) power of $q$ in the generic degree of $E$ \cite[\S 4.1.1]{Lu1}.

By \cite[Theorem 5.4 and Corollary 6.3 (b)]{Lu2},
Lusztig's $a$-function is a sorting function.
The corresponding filtrations on projective modules have been considered
in \cite{GeRou}.

Write $E<E'$ for the order on $\CF$
arising from $<_{KL}$  via the canonical
bijection between $\CC$ and $\CF$.
The following Lemma is a classical result~:
\begin{lemma}
Let $E,E'\in\CF$. If $E<E'$, then
$a_E>a_{E'}$ and $A_E >A_{E'}$.
\end{lemma}

\begin{proof}
By \cite[Remark 3.3(a)]{KaLu}, we have $v\le_{LR}w$ if and only if
$w_0w\le_{LR}w_0v$, where $w_0$ is the element of maximal length.
Left multiplication by $w_0$ induces a automorphism of $\CC$.
The corresponding automorphism of $\CF$ is tensor product by $\det$
\cite[Lemma 5.14]{Lu1}.
It follows that $E<E'$ if and only if $E'\otimes\det<E\otimes\det$
(cf also \cite[Proposition 2.25]{BaVo}).

We have $A_E=N-a_{E\otimes\det}$, where $N$ is the number of positive
roots of $W$ \cite[5.11.5]{Lu1}.

The Lemma is now a consequence of the fact that $E<E'\Rightarrow
a_E>a_{E})$.
\end{proof}

We deduce there is another sorting function~:

\begin{prop}
The function $a_E+A_E$ is a sorting function.
\end{prop}

\pagebreak[3]
\pagebreak[3]\subsection{Standard modules for the Hecke algebra via $\kz$-functor}

\subsubsection{}
We consider the setting of \S \ref{sectionreminder} with
$k_{H,1}=k_1$ independent of $H$. According to \cite[\S 4.21 and Proposition 4.1]{BrMi}
we have
$$c_E=k_1(a_E+A_E).$$

We can finally identify the images of the standard modules $\Delta(E)$ of
$\CO$ via ${\kz}$~:

\begin{thm}\label{isospecht}
Assume $k_{H,1}$ is a positive real number independent of $H$.
Let $E\in\Irr(\BC W)$. Then, ${\kz}(\Delta(E))\iso S(E)$.
\end{thm}

\begin{proof}
We prove the result for $R$ local complete as in
\S \ref{definitionHecke} instead of $\BC$. The $\Delta$-filtration of
projective objects of $\CO_R$ becomes, via ${\kz}_R$, the filtration 
of Proposition \ref{filtrationprojHecke} for the sorting function
$f(E)=a_E+A_E$
and the associated quotients are direct sums of standard $\CH_R$-modules.
It follows that the modules ${\kz}_R(\Delta(R\otimes E))$ for $E\in\Irr(W)$
coincide with the standard $\CH_R$-modules.
Since ${\kz}_K(\Delta(K\otimes E))\iso K\otimes S(E)$
(cf the remark before Theorem \ref{factorHecke})), we deduce the Theorem.
\end{proof}

\begin{rem} If the number $k_{H,1}$ (which is independent of $H$) is
  non-real,
then  the category $\CO$ and the algebra $\CH$ are both
semi-simple, hence it is still true
that  ${\kz}(\Delta(E))\simeq S(E)$. If the number is non-positive real,
then a similar approach shows that ${\kz}(\Delta(E))\simeq S(E)^*$.
\end{rem}

\begin{cor}
Assume $k_{H,1}$ is a positive real number and $W$ has type $A_n$.
Then, ${\kz}(\Delta(E))$ is isomorphic to the Specht module corresponding
to $E$.
\end{cor}

\begin{proof}
The result is a consequence of \cite{Na,GaMc}, where it is proven that
the module $S(E)$ is a Specht module.
Alternatively,
any projective $\CH$-module is known to have a
filtration by Specht modules such that the  order of terms in the
filtration
is compatible with the dominance order on partitions.
The claim of the Corollary can  be easily deduced from this
by
comparing  with the  order relation on two-sided cells.
\end{proof}

\begin{cor}\label{mult}
If $\kz(L(E))\neq 0$, then $\kz(P(E))$  is a projective $\CH$-module
and, for any $F\in \Irr(W)$, we have $[S(F):\kz(L(E))]=[P(E):\Delta(F)]$.
\end{cor}

\begin{proof} This is an immediate consequence of Theorem
\ref{isospecht}, the reciprocity
 formula in \S \ref{reciproc}, and Proposition \ref{classdeltanabla}.
\end{proof}

\smallskip

{\footnotesize
V. Ginzburg: Department of Mathematics, University of Chicago,
Chicago, IL
60637, USA;\\
\hphantom{x}\quad\, {\tt ginzburg@math.uchicago.edu}
\smallskip

N. Guay: Department of Mathematics, University of Chicago,
Chicago, IL
60637, USA;\\
\hphantom{x}\quad\, {\tt nguay@math.uchicago.edu}
\smallskip

E. Opdam: Korteweg de Vries Institute for Mathematics, University of Amsterdam,\\
\hphantom{x}\quad\,\,Plantage Muidergracht 24, 1018TV Amsterdam, The Netherlands;
\; {\tt opdam@science.uva.nl}
\smallskip

R. Rouquier: UFR de Math\'ematiques et Institut de Math\'ematiques de
Jussieu 
(CNRS UMR 7586), \\
\hphantom{x}\quad\,\,
Universit\'e Paris VII, 2 place Jussieu, 75251 Paris Cedex 05, France;
\; {\tt rouquier@math.jussieu.fr}
}

\end{document}